\documentclass[12pt]{paper}

\author{Benjamin J. Wilson}
\title{Highest-Weight Theory for Truncated Current Lie Algebras}

\usepackage{graphicx}
\usepackage{lscape}

\theoremstyle{definition}

\newtheorem*{theoremnonum}{Theorem}

\newtheorem{theorem}[equation]{Theorem}
\newtheorem{propn}[equation]{Proposition}
\newtheorem{lemma}[equation]{Lemma}
\newtheorem{corollary}[equation]{Corollary}

\newtheorem{example}[equation]{Example}

\newcommand{\K}{\Bbbk}
\newcommand{\RealNumbers}{\mathbb{R}}
\newcommand{\ComplexNumbers}{\mathbb{C}}
\newcommand{\Z}{\mathbb{Z}}
\newcommand{\N}{\mathbb{N}}
\newcommand{\Id}{\mathrm{id}}
\newcommand{\Factorial}[1]{#1 !}
\newcommand{\NewTerm}[1]{\textit{#1}}
\newcommand{\Indet}[1]{\mathrm{#1}} 
\newcommand{\Note}[1]{\text{\small{\ (#1)}}}
\newcommand{\ENote}[1]{\text{\small{\quad (#1)}}}
\newcommand{\SFrac}[2]{\textstyle{\frac{#1}{#2}}}
\newcommand{\Span}{\text{span}}

\newcommand{\Dim}{\text{dim} \hspace{0.15em}}
\newcommand{\LieBrac}[2]{\lbrack \hspace{0.15em} #1 , #2 \hspace{0.15em}\rbrack} 

\newcommand{\Ptns}{\mathcal P} 
\newcommand{\Q}{\mathcal{Q}}
\newcommand{\QPlus}{\mathcal{Q}_{+}}

\newcommand{\LinearDual}[1]{#1^{*}}
\newcommand{\g}{\mathfrak{g}}
\newcommand{\Heis}{\mathfrak{a}}
\newcommand{\h}{\mathfrak{h}}
\newcommand{\hZero}{\h_{0}}
\newcommand{\Roots}{\Delta}
\newcommand{\PositiveRoots}{\Roots_{+}}
\newcommand{\NegativeRoots}{\Roots_{-}}
\newcommand{\RootSpace}[2]{{#1}^{#2}}
\newcommand{\Implies}{\ \Rightarrow \ }
\newcommand{\DegreeFn}[1]{\text{deg}_{#1}}
\newcommand{\Degree}[2]{\text{deg}_{#1}{#2}}
\newcommand{\RootFn}{\Delta}
\newcommand{\Invol}{\omega}
\newcommand{\C}{\mathcal{C}}
\newcommand{\CHat}{\hat{\C}}
\newcommand{\gHat}{\hat{\g}}
\newcommand{\hHat}{\hat{\h}}
\newcommand{\Weight}{\Delta}
\newcommand{\Length}[1]{\vert #1 \vert}
\newcommand{\Star}[1]{{{#1}^{\star}}}
\newcommand{\UEA}[1]{\text{U}(#1)}
\newcommand{\SymAlg}[1]{\text{S}(#1)}
\newcommand{\TensorAlg}[1]{\text{T}(#1)}
\newcommand{\ShapForm}{\textbf{F}}
\newcommand{\ModShap}{\textbf{B}}
\newcommand{\Form}[2]{\ModShap \bigl(\y{#1}, \y{#2}\bigr) }
\newcommand{\Proj}{\textbf{q}}
\newcommand{\Det}[1]{\text{det}\hspace{0.1em} #1}
\newcommand{\Sym}[1]{\text{Sym}(#1)}
\newcommand{\MatrixSpace}{\mathcal L}
\newcommand{\Verma}[1]{\textrm{M}(#1)}
\newcommand{\VermaGen}[1]{\textrm{v}_{#1}}
\newcommand{\Radical}{\text{Rad}\ }
\newcommand{\x}[1]{ \mathrm{x}(#1) }
\newcommand{\y}[1]{ \mathrm{y}(#1) }
\newcommand{\ie}{i.~e.~}
\newcommand{\cf}{cf.~}

\newcommand{\WeightSpace}[2]{{#1}^{#2}}
\newcommand{\Nilp}{\mathrm{N}} 
\newcommand{\hElementMap}{\textbf{h}}
\newcommand{\hElement}[1]{\hElementMap(#1)} 
\newcommand{\PairForm}[3]{\left(#1, #2 \right)_{#3}}
\newcommand{\VertForm}[2]{(\hspace{0.2em}#1 \hspace{0.2em} \vert \hspace{0.2em}#2 \hspace{0.2em})}
\newcommand{\SymPtnPair}[2]{\textrm{S}\big(#1, #2\big)} 
\newcommand{\SymSpace}[1]{\mathrm{A}_{#1}}
\newcommand{\SymSpaceBasis}[1]{\bar{\mathrm{x}}(#1)}
\newcommand{\PtnSpanForm}[1]{\mathrm{J}_{#1}}
\newcommand{\TensorProduct}[2]{ {\textstyle #1 \bigotimes #2 }}
\newcommand{\TensorPower}[2]{\text{T}^{#1}(#2)}
\newcommand{\SymPower}[2]{\text{S}^{#1}(#2)}
\newcommand{\Reverse}[1]{{#1}^\dagger}
\newcommand{\OrderedMatrixSpace}{\MatrixSpace(\leqslant)}
\newcommand{\SpecialLinear}[1]{\mathrm{sl}(#1)}
\newcommand{\HT}[2]{\mathrm{h}_{{#1}, #2} }
\newcommand{\slh}[1]{\mathrm{h}(#1)}

\newcommand{\Vir}[2]{\Omega_{#1, #2}}
\newcommand{\VirL}{\mathrm{L}}
\newcommand{\VirC}{\mathrm{c}}
\newcommand{\VirRoot}{\updelta} 

\newcommand{\HeisA}{\mathrm{a}}
\newcommand{\HeisD}{\mathrm{d}}
\newcommand{\HeisH}{\hbar}
\newcommand{\HeisRoot}{\updelta}

\newcommand{\FiniteSeries}[2]{\mathrm{#1}_{#2}}
\newcommand{\Functional}[2]{\langle #1, #2 \rangle}
\newcommand{\Res}[1]{{#1}_\Nilp}
\newcommand{\GradedComp}[2]{{#1}_{#2}}
\newcommand{\Finite}[1]{\dot{#1}}
\newcommand{\AffineC}{\mathrm{c}}
\newcommand{\AffineD}{\mathrm{d}}
\newcommand{\ImagRoot}{\updelta}
\newcommand{\OtherRoot}{\uptau}
\newcommand{\FiniteRes}[1]{\widetilde{#1}}
\newcommand{\RealRoots}{\Roots^{\text{re}}}
\newcommand{\ImaginaryRoots}{\Roots^{\text{im}}}

\newcommand{\Character}[1]{\textrm{char}\hspace{0.15em} #1}
\newcommand{\Irred}[1]{\textrm{L}(#1)}
\newcommand{\Exp}[1]{\textrm{e}^{#1}}

\begin{document}
\begin{abstract}
Let $\g$ be a Lie algebra over a field $\K$ of characteristic zero, and a fix positive integer $\Nilp$.
The Lie algebra
$$ \gHat = \g \otimes_\K \K[\Indet{t}] / {\Indet{t}^{\Nilp+1} \K[\Indet{t}]}$$
is called a truncated current Lie algebra.
In this paper a highest-weight theory for $\gHat$ is developed when the underlying Lie algebra $\g$ possesses a triangular decomposition.
The principal result is the reducibility criterion for the Verma modules of $\gHat$ for a wide class of Lie algebras $\g$, including the symmetrizable Kac-Moody Lie algebras, the Heisenberg algebra, and the Virasoro algebra.
This is achieved through a study of the Shapovalov form.
\end{abstract}

\maketitle
\begin{section}{Overview}
Let $\g$ be a Lie algebra over a field $\K$ of characteristic zero, and fix a positive integer $\Nilp$.
The Lie algebra
\begin{equation}\label{AffinizationEqn}
\gHat = \g \otimes_\K \K[\Indet{t}] / {\Indet{t}^{\Nilp+1} \K[\Indet{t}]}
\end{equation}
over $\K$, with the Lie bracket
\begin{equation}\label{Intro:LieBracket}
\LieBrac{x \otimes \Indet{t}^i}{y \otimes \Indet{t}^j} = \LieBrac{x}{y} \otimes \Indet{t}^{i+j}, \qquad x,y \in \g, \quad i,j \geqslant 0,
\end{equation}
is called a \NewTerm{truncated current Lie algebra}.
The Lie algebra $\gHat$ is graded by non-negative degrees in $\Indet{t}$.
Any graded representation of $\gHat$ with more than one graded component is trivially reducible, and so we disregard the grading in $\Indet{t}$ in our representation theory.
A triangular decomposition (\cf Section \ref{LATDSection}, here $\hZero = \h$)
\begin{equation}\label{Intro:TriangularDecomposition}
\g = \g_{-} \oplus \h \oplus \g_{+}, \qquad \g_{\pm} = \oplus_{\alpha \in \PositiveRoots}{\g^{\pm \alpha}}, \qquad \PositiveRoots \subset \LinearDual{\h}, 
\end{equation}
of $\g$ naturally defines a triangular decomposition of $\gHat$,
$$ \gHat = \gHat_{-} \oplus \hHat \oplus \gHat_{+}, \qquad \gHat_{\pm} = \oplus_{\alpha \in \PositiveRoots }{\gHat^{\pm \alpha}},$$
where the subalgebra $\hHat$ and the subspaces $\RootSpace\gHat\alpha$ are defined by analogy with (\ref{AffinizationEqn}), and $\h \subset \hHat$ is the diagonal subalgebra.
Thus a weight module for $\gHat$ is a $\gHat$-module that is $\h$-diagonalizable, while a highest-weight module for $\gHat$ is a weight $\gHat$-module generated by an eigenvector for $\hHat$ that is annihilated by $\gHat_{+}$.
In this paper we construct a highest-weight theory for $\gHat$, and derive reducibility criteria for the universal objects of the theory, the Verma modules, by studying the Shapovalov form.

We assume that the decomposition (\ref{Intro:TriangularDecomposition}) of $\g$ is non-degenerately paired, \ie that for each $\alpha \in \PositiveRoots$, a non-degenerate bilinear form
$$ {\VertForm\cdot\cdot}_\alpha : \RootSpace\g\alpha \times \RootSpace\g{-\alpha} \to \K,$$
and a non-zero element $\hElement\alpha \in \h$ are given, such that
$$ \LieBrac{x}{y} = {\VertForm{x}{y}}_\alpha \hspace{0.15em} \hElement\alpha,$$
for all $x \in \RootSpace\g\alpha$ and $y \in \RootSpace\g{-\alpha}$ (\cf Section \ref{BlockSection} for a precise definition).
The symmetrizable Kac-Moody Lie algebras, the Virasoro algebra and the Heisenberg algebra all possess triangular decompositions that are non-degenerately paired.
The principal result of this paper is the following.\par
\vspace{0.5cm}
\begin{theoremnonum}
A Verma module $\Verma\Lambda$ for $\gHat$ is reducible if and only if
$$ \Functional\Lambda{\hElement\alpha \otimes \Indet{t}^\Nilp} = 0$$
for some positive root $\alpha \in \PositiveRoots$ of $\g$.
\end{theoremnonum}
\vspace{0.5cm}
Truncated current Lie algebras have previously been studied in the case where $\g$ is a semisimple finite-dimensional Lie algebra.
Takiff considered this case with $\Nilp=1$ in \cite{Takiff}, and that work was extended in \cite{RaisTauvel}, \cite{GeoffriauOne}, \cite{GeoffriauTwo} without the restriction on $\Nilp$.
As such, when $\g$ is a semisimple finite-dimensional Lie algebra, the Lie algebra $\gHat$ is known as a \NewTerm{generalized Takiff algebra}.

An \NewTerm{exp-polynomial module} for a loop algebra
$$ \tilde\g = \g \otimes \K [ \Indet{t}, \Indet{t}^{-1}]$$
(with the Lie bracket given by (\ref{Intro:LieBracket})) is a module for which the defining relations are equations involving only polynomial and exponential dependence upon the power of the indeterminate~$\Indet{t}$.
Berman, Billig and Zhao proved in \cite{BermanBillig}, \cite{BilligZhao} that any irreducible exp-polynomial module for a loop algebra must have finite-dimensional graded components.
The derivation of the character formulae is an interesting open problem.

The irreducible exp-polynomial modules for loop algebras may be realized through affinizations of irreducible highest-weight modules for truncated current Lie algebras.
Character formulae for Verma modules are relatively straightforward to derive.
Thus the reducibility criterion detailed above allows the calculation of the characters of the irreducible exp-polynomial modules in many cases. 
The calcuation of these characters is the current focus of the author's research.
\end{section}

\begin{section}{Lie Algebras with Triangular Decompositions}\label{LATDSection}
Let $\g$ be a Lie algebra over a field $\K$.
A \NewTerm{triangular decomposition} of $\g$ is specified by a pair of non-zero abelian subalgebras $\hZero \subset \h$, a pair of distinguished non-zero subalgebras $\g_{+}$, $\g_{-}$, and an anti-involution (\ie an anti-automorphism of order $2$) 
$$\Invol: \g \to \g,$$ 
such that:
\begin{enumerate}
\item $\g = \g_{-} \oplus \h \oplus \g_{+}$;
\item the subalgebra $\g_{+}$ is a non-zero weight module for $\hZero$ under the adjoint action, with weights $\PositiveRoots$ all non-zero; 
\item $\Invol|_\h = \Id_\h$ and $\Invol(\g_{+}) = \g_{-}$;
\item the semigroup with identity $\QPlus$, generated by $\PositiveRoots$ under addition, is freely generated by a finite subset
$\set{\alpha_j}_{j \in {\mathrm J}} \subset \QPlus$
consisting of linearly independent elements of $\LinearDual{\hZero}$.
\end{enumerate}
This definition is a modification of the definition of Moody and Pianzola \cite{MoodyPianzola}.
There, the set $\mathrm J$ is not required to be finite, root spaces may be infinite-dimensional, and $\hZero = \h$.
We distinguish between $\hZero$ and $\h$ in order to include Example \ref{TrunCurrentExample}.

Write $\Q = \sum_{j \in {\mathrm J}}{\Z \alpha_j}$.
Call the weights $\PositiveRoots$ of the $\hZero$-module $\g_{+}$ the \NewTerm{positive roots}, and the weight space $\RootSpace{\g}{\alpha}$ corresponding to $\alpha \in \PositiveRoots$ the \NewTerm{$\alpha$-root space}, so that $\g_{+} = \oplus_{\alpha \in \PositiveRoots}{\RootSpace{\g}{\alpha}}$.
The anti-involution ensures an analogous decomposition of $\g_{-} = \oplus_{\alpha \in \NegativeRoots}{\RootSpace{\g}{\alpha}}$, where $\NegativeRoots = - \PositiveRoots$ (the \NewTerm{negative roots}) and $\RootSpace{\g}{-\alpha} = \Invol(\RootSpace{\g}{\alpha})$ for all $\alpha \in \PositiveRoots$.
Write $\Roots = \PositiveRoots \cup \NegativeRoots$ for the \NewTerm{roots} of $\g$.
Consider $\QPlus$ to be partially ordered in the usual manner, i.e. for $\gamma, \gamma' \in \QPlus$,
$$ \gamma \leqslant_{\QPlus} \gamma' \quad \iff \quad (\gamma' - \gamma) \in \QPlus.$$
We assume that all root spaces are finite-dimensional,
and that $\PositiveRoots$ is a countable set.
For clarity, a Lie algebra with triangular decomposition may be referred to as a five-tuple $(\g, \hZero, \h, \g_{+}, \Invol)$.

\begin{example}\label{FDLAExample}
Let $\g$ be a finite-dimensional semisimple Lie algebra over $\ComplexNumbers$, with Cartan subalgebra $\h$ and root system $\Roots$.  Then
$$ \g = \h \oplus (\oplus_{\alpha \in \Roots} \RootSpace{\g}{\alpha}).$$
Let $\pi$ be a basis for $\Roots$, and let $\QPlus$ be the additive semigroup generated by $\pi$.
Write $\PositiveRoots = \Roots \cap \QPlus$, and let $\g_{+}$, $\g_{-}$ be given by
$$ \g_{+} = \oplus_{\alpha \in \PositiveRoots} \RootSpace{\g}{\alpha}, \qquad \g_{-} = \oplus_{\alpha \in \NegativeRoots} \RootSpace{\g}{\alpha},$$
where $\NegativeRoots = - \PositiveRoots$.
Then $\g_{+}$ is a weight-module for $\hZero = \h$ with weights $\PositiveRoots$, and 
$$ \g = \g_{-} \oplus \h \oplus \g_{+}.$$
All root spaces are one-dimensional.
For any $\alpha \in \PositiveRoots$, choose non-zero elements
$$ \x\alpha \in \RootSpace{\g}{\alpha}, \qquad \y\alpha \in \RootSpace{\g}{\alpha}.$$
An anti-involution $\Invol$ on $\g$ is defined by extension of
$$ \Invol|_{\h} = \Id_{\h}, \qquad \Invol(\x\alpha) = \y\alpha, \quad \Invol(\y\alpha) = \x\alpha, \quad \alpha \in \pi.$$
Thus $(\g, \g_{+}, \h, \h, \Invol)$ is a Lie algebra with triangular decomposition.
The semisimple finite-dimensional Lie algebras over $\ComplexNumbers$ are parameterized by Euclidean root systems, or equivalently by the Cartan matrices.
The Serre relations permit the construction of any such Lie algebra from its Cartan matrix, and this construction works over an arbitrary field $\K$ of characteristic zero.
The preceding assertions hold also for the Lie algebras over $\K$ constructed in this manner.
Here and throughout, \NewTerm{semisimple finite-dimensional Lie algebra} means a Lie algebra over $\K$ defined by a Cartan matrix and the Serre relations.
\end{example}

\begin{example}\label{SL3Example}
It shall be convenient to consider the following particular case of Example \ref{FDLAExample} is greater detail.
Let $\g$ denote $\SpecialLinear{3}$, the finite-dimensional semisimple Lie algebra over $\K$ with root system $\FiniteSeries{A}{2}$.
Denote by $\upalpha_1$, $\upalpha_2$ the simple roots, by 
$$\x{\upalpha_1}, \ \x{\upalpha_2}, \ \y{\upalpha_1}, \ \y{\upalpha_2}, \ \slh{\upalpha_1}, \ \slh{\upalpha_2} $$
the Chevalley generators,
and by $\h = \K \slh{\upalpha_1} \oplus \K \slh{\upalpha_2}$ the Cartan subalgebra, so that
$$ \upalpha_1 ( \slh{\upalpha_1} ) = \upalpha_2 ( \slh{\upalpha_2} ) = 2, \quad \upalpha_1 ( \slh{\upalpha_2} ) = \upalpha_2 ( \slh{\upalpha_1} ) = -1.$$ 
Then the root system is defined by $\Roots = \PositiveRoots \cup \NegativeRoots$, where $\PositiveRoots = \set{ \upalpha_1, \upalpha_2, \upalpha_1 + \upalpha_2 }$ and $\NegativeRoots = - \PositiveRoots$.
Write
$$ \x{\upalpha_1 + \upalpha_2} = \LieBrac{\x{\upalpha_1}}{\x{\upalpha_2}}, \quad \y{\upalpha_1 + \upalpha_2} = \LieBrac{\y{\upalpha_2}}{\y{\upalpha_1}}, \quad \slh{\upalpha_1 + \upalpha_2} = \slh{\upalpha_1} + \slh{\upalpha_2}.$$
Then for each $\alpha \in \PositiveRoots$, the elements $\x\alpha, \y\alpha, \slh\alpha$ span a subalgebra of $\g$ isomorphic to $\SpecialLinear{2}$.
The anti-involution $\Invol$ fixes $\h$ point-wise, and interchanges $\x\alpha$ with $\y\alpha$ for every $\alpha \in \PositiveRoots$.
Write 
$$\g^\alpha = \K \x\alpha, \quad \g^{-\alpha} = \K \y\alpha, \quad \alpha \in \PositiveRoots, $$
and $\g_{\pm} = \oplus_{\alpha \in \PositiveRoots} \g^{\pm \alpha}$.  Then $\g_{+}$ is a weight-module for $\hZero = \h$ with weights $\PositiveRoots$, and 
$$ \g = \g_{-} \oplus \h \oplus \g_{+}.$$
The semigroup $\QPlus$ is generated by $\pi = \set{\upalpha_1, \upalpha_2}$.
Note that the $\slh{\alpha}$ defined here are only proportional to the elements $\hElement{\alpha}$ defined later on.
\end{example}

\begin{example}\label{KMLAExample}
Let $\g$ be the Kac-Moody Lie algebra over $\K$ associated to an $n \times n$ generalized Cartan matrix (we paraphrase \cite{KacBook}).
Let $\h$ denote the Cartan subalgebra, and $\Roots$ the root system.
Then
$$ \g = \h \oplus (\oplus_{\alpha \in \Roots} \RootSpace{\g}{\alpha}),$$
and all root spaces are finite-dimensional.
The collection $\Pi$ of simple roots is a linearly-independent subset of the finite-dimensional space $\LinearDual\h$.
Let $\QPlus$ denote the additive semigroup generated by $\Pi$, let $\PositiveRoots = \Roots \cap \QPlus$, and write
$$ \g_{+} = \oplus_{\alpha \in \PositiveRoots} \RootSpace{\g}{\alpha}, \qquad \g_{-} = \oplus_{\alpha \in \NegativeRoots} \RootSpace{\g}{\alpha},$$
where $\NegativeRoots = - \PositiveRoots$.
Then $\g_{+}$ is a weight-module for $\hZero = \h$ with weights $\PositiveRoots$, and 
$$ \g = \g_{-} \oplus \h \oplus \g_{+}.$$
If $e_i, f_i$, $1 \leqslant i \leqslant n$, denote the Chevalley generators of $\g$, then $\g_{+}$ and $\g_{-}$ are the subalgebras generated by the $e_i$ and by the $f_i$, respectively.
An anti-involution $\Invol$ of $\g$ is defined by extension of
$$ \Invol|_{\h} = \Id_{\h}, \qquad \Invol(e_i) = f_i, \quad \Invol(f_i) = e_i, \quad 1 \leqslant i \leqslant n,$$
(this $\Invol$ differs from the $\omega$ of \cite{KacBook}).
Thus $(\g, \g_{+}, \h, \h, \Invol)$ is a Lie algebra with triangular decomposition.
\end{example}

\begin{example}\label{VirasoroExample}
Let $\g$ denote the $\K$-vector space with basis the symbols 
$$\set{\VirL_m | m \in \Z} \cup \set{ \VirC},$$
endowed with the Lie bracket given by
$$ \LieBrac{\VirC}{\g}=0, \quad \LieBrac{\VirL_m}{\VirL_n} = (m-n) \VirL_{m+n} + \delta_{m,-n} \psi(m) \VirC, \quad m,n \in \Z,$$
where $\psi : \Z \to \K$ is any function satisfying $\psi(-m) = - \psi(m)$ for $m \in \Z$, and 
\begin{equation*}
\psi (m+n) = \SFrac{2m+n}{n-m}\psi(n) + \SFrac{m + 2n}{m-n}\psi(m), \qquad m,n \in \Z, \quad m \ne n.
\end{equation*}
If $\psi = 0$, then the symbols $\VirL_m$ span a copy of the Witt algebra.
The Virasoro algebra is the only non-split one-dimensional central extension of the Witt algebra, up to isomorphism \cite{KacRaina}, and is typically defined with $\psi(m) = \SFrac{m^3 -m}{12}$.
Let 
$$\g_{\pm} = \oplus_{m > 0}{\K \VirL_{\pm m}}, \quad \hZero = \h = {\K \VirL_0} \oplus {\K \VirC},$$ 
and let $\VirRoot \in \LinearDual{\h}$ be given by
$$ \VirRoot (\VirL_0) = -1, \quad \VirRoot (\VirC) = 0.$$
Then $\g = \g_{-} \oplus \h \oplus \g_{+}$, and $\g_{+}$ is a weight module for $\hZero = \h$, with weights 
$$\PositiveRoots = \set{m \VirRoot | m > 0}.$$
The semigroup $\QPlus$ is generated by $\VirRoot$.
An anti-involution $\Invol$ is given by 
$$\Invol (\VirC) = \VirC, \quad \Invol(\VirL_m) = \VirL_{-m}, \quad m \in \Z,$$
and in this notation $\g$ is a Lie algebra with triangular decomposition.
\end{example}

\begin{example}\label{HeisenbergExampleOne}
Let $\Heis$ denote the $\K$-vector space with basis the symbols 
$$\set{\HeisA_m | m \in \Z} \cup \set{  \HeisH, \ \HeisD},$$ 
endowed with the Lie bracket given by
$$ \LieBrac{\HeisA_m}{\HeisA_n} = m \delta_{m,-n} \HeisH, \quad \LieBrac{\HeisH}{\Heis} = 0, \quad \LieBrac{\HeisD}{\HeisA_m} = m \HeisA_m, \quad m,n \in \Z.$$
The Lie algebra $\Heis$ is called the extended Heisenberg or oscillator algebra.
Let
$$ \Heis_{\pm} = \oplus_{m>0}{\K \HeisA_{\pm m}}, \quad \h = \K \HeisA_0 \oplus \K \HeisH \oplus \K \HeisD,$$
and let $\HeisRoot \in \LinearDual{\h}$ be given by
$$ \HeisRoot (\HeisA_0) = \HeisRoot (\HeisH) = 0, \quad \HeisRoot(\HeisD) = 1.$$
Then $\Heis = \Heis_{-} \oplus \h \oplus \Heis_{+}$, and $\Heis_{+}$ is a weight module for $\hZero = \h$, with weights 
$$\PositiveRoots = \set{m \HeisRoot | m > 0}.$$
The semigroup $\QPlus$ is generated by $\HeisRoot$.
An anti-involution $\Invol$ is given by 
$$\Invol(\HeisH) = \HeisH, \quad \Invol(\HeisD) = \HeisD, \quad \Invol(\HeisA_m) = \HeisA_{-m}, \quad m \in \Z,$$
and in this notation $\Heis$ is a Lie algebra with triangular decomposition.
\end{example}

\begin{example}\label{TrunCurrentExample}
Let $\g$ be a $\K$-Lie algebra with triangular decomposition, denoted as above, and let $R$ be a commutative, associative $\K$-algebra with $1$ (e.g. $R = \K[\Indet{t}] / \Indet{t}^{\Nilp+1} \K[\Indet{t}]$, $\Nilp>0$).
Write $\gHat = \g \otimes_{\K} R$, and similarly for the subalgebras of $\g$.
Then $\gHat$ is a $\K$-Lie algebra with Lie bracket
$$ \LieBrac{x \otimes r}{y \otimes s} = \LieBrac{x}{y}\otimes{rs}, \quad x,y \in \g, \ r,s \in R,$$
and contains $\g$ as a subalgebra via $x \mapsto x \otimes 1$.
Moreover, $ \gHat = \gHat_{-} \oplus \hHat \oplus \gHat_{+}$, and $\hZero \subset \hHat$ are non-zero abelian subalgebras of $\gHat$.
The subalgebra $\gHat_{+}$ is a weight module for $\hZero$ with weights coincident with the weights $\PositiveRoots$ of the $\hZero$-module $\g_{+}$, and $(\gHat_{+})^\alpha = \widehat{ (\g_{+}^\alpha)}$.
So $\g$ and $\gHat$ share the same roots $\Roots$ and root lattices $\Q$, $\QPlus$.
The anti-involution $\Invol$ of $\gHat$ is given by $R$-linear extension
$$ \Invol : x \otimes r \mapsto \Invol(x) \otimes r, \quad x \in \g, \ r \in R,$$
and fixes $\hHat$ point-wise.
Thus $(\gHat, \hZero, \hHat, \gHat_{+}, \Invol)$ is a $\K$-Lie algebra with triangular decomposition.
\end{example}

\end{section}

\begin{section}{Highest-Weight Theory for Lie Algebras with Triangular Decomposition}\label{HWTSection}
Throughout this section, let $(\g, \hZero, \h, \g_{+}, \Invol)$ denote a Lie algebra with triangular decomposition.
The universal highest-weight modules of $\g$, called \NewTerm{Verma modules}, exist and possess the usual properties.
An extensive treatment of Verma modules and the Shapovalov form can be found in \cite{MoodyPianzola}; we present only the definitions and the most important properties.
\begin{subsection}{Highest-weight modules}
A $\g$-module $M$ is \NewTerm{weight} if the action of $\hZero$ on $M$ is diagonalizable, \ie
\begin{equation}\label{WSDecomp}
M = \oplus_{\chi \in \LinearDual{\hZero}}{\WeightSpace{M}{\chi}}, \qquad h|_{\WeightSpace{M}{\chi}} = \chi (h) \text{\ \ for all\ \ } h \in \hZero, \chi \in \LinearDual{\hZero}.
\end{equation}
The decomposition (\ref{WSDecomp}) is called the \NewTerm{weight-space decomposition} of $M$; the components $\WeightSpace{M}{\chi}$ are called \NewTerm{weight spaces}.
The \NewTerm{weight lattice} of a weight module $M$ is the set
$$ \set{ \chi \in \LinearDual{\hZero} | \WeightSpace{M}{\chi} \ne 0} \subset \LinearDual{\hZero}.$$
For any $\chi \in \LinearDual\hZero$, an element $v \in \WeightSpace{M}{\chi}$ is a \NewTerm{primitive vector} of $M$ if the submodule $\UEA\g \cdot v \subset M$ is proper. 
Clearly $M$ is reducible if and only if $M$ has a non-zero primitive vector.
A non-zero vector $v \in M$ is a \NewTerm{highest-weight vector} if 
\begin{enumerate}
\item $\g_{+} \cdot v = 0$;
\item there exists $\Lambda \in \LinearDual{\h}$ such that $ h \cdot v = \Lambda (h) v$, for all $h \in \h$. 
\end{enumerate}
The unique functional $\Lambda \in \LinearDual{\h}$ is called the \NewTerm{highest weight} of the highest-weight vector $v$.
A weight $\g$-module $M$ is called \NewTerm{highest weight} (of highest weight $\Lambda$) if there exists a highest-weight vector $v \in M$ (of highest weight $\Lambda$) that generates it.

\begin{propn}\label{HWMPropn}
Suppose that $M$ is a highest-weight $\g$-module, generated by a highest-weight vector $v \in M$ of highest weight $\Lambda \in \LinearDual\h$. 
Then
\begin{enumerate}
\item the weight lattice of $M$ is contained in 
$\Lambda|_{\hZero} - \QPlus$; 
\item $\WeightSpace{M}{\Lambda|_{\hZero}} = \K v$, and all weight spaces of $M$ are finite-dimensional;
\item $M$ is indecomposable, and has a unique maximal submodule;
\item if $u \in M$ is a highest-weight vector of highest-weight $\Lambda' \in \LinearDual{\h}$, and $u$ generates $M$, then $\Lambda' = \Lambda$ and $u$ is proportional to $v$.
\end{enumerate}
\end{propn}
Let $\Lambda \in \LinearDual{\h}$, and
consider the one-dimensional vector space $\K \VermaGen\Lambda$ as an $(\h \oplus \g_{+})$-module via
$$ \g_{+} \cdot \VermaGen\Lambda = 0; \quad h \cdot \VermaGen\Lambda = \Lambda (h) \VermaGen\Lambda, \ h \in \h.$$
The induced module
$$ \Verma{\Lambda} = \UEA{\g} \otimes_{\UEA{\h \oplus \g_{+}}} \K \VermaGen\Lambda$$
is called the \NewTerm{Verma module} of highest-weight $\Lambda$.
\begin{propn}\label{VermaPropertiesPropn}
For any $\Lambda \in \LinearDual{\h}$, 
\begin{enumerate}
\item\label{VermaPropertiesPropn1} Up to scalar multiplication, there is a unique epimorphism from $\Verma{\Lambda}$ to any highest-weight module of highest-weight $\Lambda$, i.e. $\Verma{\Lambda}$ is the universal highest-weight module of highest-weight $\Lambda$;
\item\label{VermaPropertiesPropn2} $\Verma{\Lambda}$ is a free rank one $\UEA{\g_{-}}$-module.
\end{enumerate}
\end{propn}
\end{subsection}
\begin{subsection}{The Shapovalov Form}
The Shapovalov form is a contragredient symmetric bilinear form on $\UEA{\g}$ with values in $\UEA{\h} = \SymAlg{\h}$.
The evaluation of the Shapovalov form at $\Lambda \in \LinearDual\h$ is a $\K$-valued bilinear form, and is degenerate if and only if the Verma module $\Verma\Lambda$ is reducible.
By the Leibniz rule, $\UEA{\g}$ is a weight $\g$-module, with weight-space decomposition
$$ \UEA{\g} = \oplus_{\gamma \in \Q}\WeightSpace{\UEA{\g}{\gamma}}.$$
The anti-involution $\Invol$ of $\g$ extends uniquely to an anti-involution of $\UEA{\g}$ (denoted identically), and is such that 
$$\Invol : \UEA{\g}^{\gamma} \to \UEA{\g}^{- \gamma}, \quad \gamma \in \Q.$$
It follows from the Poincar\'e-Birkhoff-Witt (PBW) theorem that $\UEA{\g}$ may be decomposed
$$ \UEA{\g} = \UEA{\h} \oplus \{ \g_{-} \UEA{\g} + \UEA{\g} \g_{+} \}$$
as a direct sum of vector spaces.
Further, both summands are two-sided $\UEA{\h}$-modules preserved by $\Invol$.
Let $\Proj : \UEA{\g} \to \UEA{\h}$ denote the projection onto the first summand parallel to the second;
the restriction $\Proj |_{\UEA{\g}^0}$ is an algebra homomorphism.
Define
$$ \ShapForm : \UEA{\g} \times \UEA{\g} \to \UEA{\h} \quad \text{via} \quad \ShapForm (x,y) = \Proj (\Invol(x) y), \quad x,y \in \UEA{\g}.$$
The bilinear form $\ShapForm$ is called the \NewTerm{Shapovalov form}; we consider its restriction
$$ \ShapForm : \UEA{\g_{-}} \times \UEA{\g_{-}} \to \UEA{\h}.$$
Distinct $\hZero$-weight spaces of $\UEA{\g_{-}}$ are orthogonal with respect to $\ShapForm$, and so the study of $\ShapForm$ on $\UEA{\g_{-}}$ reduces to the study of the restrictions
$$ \ShapForm_{\chi} : \UEA{\g_{-}}^{-\chi} \times \UEA{\g_{-}}^{-\chi} \to \UEA{\h}, \quad \chi \in \QPlus. $$
Any $\Lambda \in \LinearDual{\h}$ extends uniquely to a map $\UEA{\h} \to \K$; write $\ShapForm_{\chi}(\Lambda)$ for the composition of $\ShapForm_\chi$ with this extension, and write $\Radical \ShapForm_{\chi} (\Lambda)$ for its radical.
The importance of the Shapovalov form stems from the following fact.
\begin{propn}\label{RadicalCoincidesWithMaxSubmodule}
Let $\chi \in \QPlus$, $\Lambda \in \LinearDual{\h}$.
Then $\Radical \ShapForm_{\chi} (\Lambda) \subset \Verma{\Lambda}^{\Lambda|_{\hZero} - \chi}$ is the $\Lambda|_{\hZero} - \chi$ weight space of the maximal submodule of the Verma module $\Verma{\Lambda}$.
\end{propn}
In particular, a Verma module $\Verma{\Lambda}$ is irreducible if and only if the forms $\ShapForm_{\chi}(\Lambda)$ are non-degenerate for every $\chi \in \QPlus$.
Thus an understanding of the forms $\ShapForm_\chi$, $\chi \in \QPlus$, is an understanding of the irreducibility criterion of the Verma modules of the highest-weight theory.
\end{subsection}
\begin{subsection}{Partitions and the Poincar\'e-Birkhoff-Witt Monomials}\label{PartitionsSubsection}
Let $\C$ be a set parameterizing a root-basis (i.e. an $\hZero$-weight basis) of $\g_{+}$, via
$$ \C \ni \quad \gamma \quad \leftrightarrow \quad \x{\gamma} \quad \in \g_{+}.$$
Define $\RootFn : \C \to \PositiveRoots$ by declaring 
$ \x{\gamma} \in \g_{+}^{\RootFn(\gamma)}$,
for all $\gamma \in \C$.
A \NewTerm{partition} is a finite multiset with elements from $\C$; write $\Ptns$ for the set of all partitions.
Set notation is used for multisets throughout.
In particular,  the \NewTerm{length} $\Length\lambda$ of a partition $\lambda \in \Ptns$ is the number of elements of $\lambda$, counting all repetition.
Fix some ordering of the basis $\set{\x{\gamma} | \gamma \in \C}$ of $\g_{+}$; for any $\lambda \in \Ptns$, let
\begin{equation}\label{PBWMonomialEqn}
\x\lambda = \x{\lambda_1} \cdots \x{\lambda_k} \quad \in \UEA{\g_{+}}
\end{equation}
where $k = \Length\lambda$ and $(\lambda_i)_{1 \leqslant i \leqslant k}$  is an enumeration of the entries of $\lambda$ such that (\ref{PBWMonomialEqn}) is a PBW monomial with respect to the basis ordering.
For any partition $\lambda \in \Ptns$, write $\y{\lambda} = \Invol(\x{\lambda})$.
By the PBW Theorem, the spaces $\UEA{\g_{+}}$, $\UEA{\g_{-}}$ have bases
$$ \set{\x\lambda | \lambda \in \Ptns}, \qquad \set{\y\lambda | \lambda \in \Ptns},$$
respectively.
For any partition $\lambda \in \Ptns$ and positive root $\alpha \in \PositiveRoots$, write 
$$
\RootFn (\lambda) = \sum_{\gamma \in \lambda} \RootFn (\gamma); \qquad \lambda^\alpha = \set{\gamma \in \lambda | \RootFn (\gamma) = \alpha}.
$$
\end{subsection}

\begin{subsection}{Shapovalov's Lemma}
The proof of the following useful lemma is elementary.
\begin{lemma}\label{ShuffleLemma}
Suppose that $\lambda \in \Ptns$, that $\Length{\lambda} = r$, that $(\lambda_i)_{1 \leqslant i \leqslant r}$ is an enumeration of $\lambda$ and that $\tau \in \Sym{r}$.  Then
$$ \x{\lambda_1} \cdots \x{\lambda_r} = \x{\lambda_{\tau(1)}} \cdots \x{\lambda_{\tau(r)}} + R $$
where $R$ is a linear combination of terms $\x{\phi_1} \cdots \x{\phi_s}$ where $\phi_i \in \C$ for $1 \leqslant i \leqslant s$ and $s < r$.
\end{lemma}

The following Lemma is due to Shapovalov \cite{Shapovalov}.
Our proof follows that of an analogous statement in \cite{MoodyPianzola}.

\begin{lemma}\label{ShapLemma}
Let $(\g, \hZero, \h, \g_{+}, \Invol)$ be a Lie algebra with triangular decomposition.
Suppose that $\lambda, \mu \in \Ptns$, that $\Length{\lambda} = r$ and $\Length{\mu} = s$, and that $(\lambda_i)_{1 \leqslant i \leqslant r}$ and $(\mu_i)_{1 \leqslant i \leqslant s}$ are arbitrary enumerations of $\lambda$ and $\mu$, respectively.
Let 
$$Z = \x{\lambda_r} \cdots \x{\lambda_1} \y{\mu_1} \cdots \y{\mu_s}.$$
Then
\begin{enumerate}
\item \label{ShapLemmaOne} $\Degree{\h}{\Proj (Z)} \leqslant r,s$;
\item \label{ShapLemmaTwo} if $r=s$, but $\Length{\lambda^\alpha} \ne \Length{\mu^\alpha}$ for some $\alpha \in \PositiveRoots$, then $\Degree{\h}{\Proj (Z)} < r = s$;
\item \label{ShapLemmaThree} if $r=s$ and $\Length{\lambda^\alpha} = \Length{\mu^\alpha} =: m_\alpha$ for all $\alpha \in \PositiveRoots$, then the degree $r=s$ term of $\Proj(Z)$ is 
$$ \prod_{\alpha \in \PositiveRoots}  \sum_{\tau \in \Sym{m_\alpha}} \prod_{1 \leqslant j \leqslant m_\alpha} 
	\LieBrac{ \x{\lambda^\alpha_{\tau (j)}} }{ \y{\mu^\alpha_j} },$$
where for each $\alpha \in \PositiveRoots$, $(\lambda^\alpha_j)_{1 \leqslant j \leqslant m_\alpha}$, $(\mu^\alpha_j)_{1 \leqslant j \leqslant m_\alpha}$ are any fixed enumerations of $\lambda^\alpha$ and $\mu^\alpha$ respectively.
\end{enumerate}
\end{lemma}
\begin{proof}
The proof is by induction on $\Length{\lambda} + \Length{\mu}$.
It is straightforward to show that all three parts hold whenever $\Length\lambda= 0$ or $\Length\mu = 0$.
Suppose then that all three parts hold for all $\lambda', \mu' \in \Ptns$ such that $\Length{\lambda'} + \Length{\mu'} < \Length{\lambda} + \Length{\mu}$.
Let $\varsigma \in \Sym{r}, \tau \in \Sym{s}$, and write $$Z' = \x{\lambda_{\varsigma(r)}} \cdots \x{\lambda_{\varsigma(1)}} \y{\mu_{\tau (1)}} \cdots \y{\mu_{\tau (s)}}.$$
By Lemma \ref{ShuffleLemma}, $Z = Z' + R$, where $R$ is a linear combination of terms $$\x{\phi_{r'}} \cdots \x{\phi_{1}} \y{\psi_1} \cdots \y{\psi_{s'}}$$ with $r' < r$ or $s' < s$.
Therefore, by inductive hypothesis, the Lemma holds for arbitrary enumerations of $\lambda$ and $\mu$, if it holds for any particular pair of enumerations.
Consider $\PositiveRoots$ to carry some linearization of its usual partial order, and choose any enumerations of $\lambda$, $\mu$ such that
$$ \RootFn(\lambda_1) \leqslant \cdots \leqslant \RootFn(\lambda_r) \ \ \text{and} \ \ \RootFn(\mu_1) \leqslant \cdots \leqslant \RootFn(\mu_s).$$
Moreover, as 
\begin{eqnarray*}
\Proj (\x{\lambda_r} \cdots \x{\lambda_1} \y{\mu_1} \cdots \y{\mu_s}) &=& \Proj ( \hspace{0.15em} \Invol(\x{\lambda_r} \cdots \x{\lambda_1} \y{\mu_1} \cdots \y{\mu_s}) \hspace{0.15em} ) \\
	&=& \Proj (\x{\mu_s} \cdots \x{\mu_1} \y{\lambda_1} \cdots \y{\lambda_r}),
\end{eqnarray*}
it may supposed without loss of generality that $\RootFn(\mu_1) \leqslant \RootFn(\lambda_1)$.
Now
\begin{eqnarray*}
\Proj (Z) &=& \Proj (\x{\lambda_r} \cdots \x{\lambda_1} \y{\mu_1} \cdots \y{\mu_s}) \\
	&=& \Proj ( \LieBrac{\x{\lambda_r} \cdots \x{\lambda_1}}{\y{\mu_1}} \y{\mu_2} \cdots \y{\mu_s}) \\
	&=& \sum_{i=1}^{r} \Proj (A_i),
\end{eqnarray*}
where, by the Leibniz rule,
$$
A_i = \x{\lambda_r} \cdots \x{\lambda_{i+1}} \LieBrac{\x{\lambda_i}}{\y{\mu_1}} \x{\lambda_{i-1}} \cdots \x{\lambda_1} \y{\mu_2} \cdots \y{\mu_s},$$
for $1 \leqslant i \leqslant r$.
Let $0 \leqslant k \leqslant r$ be maximal such that $\RootFn(\lambda_i) = \RootFn(\mu_1)$ for all $1 \leqslant i \leqslant k$; the terms $A_i$ with $1 \leqslant i \leqslant k$ and $k < i \leqslant r$ are to be considered separately.
If $1 \leqslant i \leqslant k$, then $\LieBrac{\x{\lambda_i}}{\y{\mu_1}} \in \h$.
Therefore, by the Leibniz rule,
$$ A_i = Z_i \LieBrac{\x{\lambda_i}}{\y{\mu_1}} + R_i,$$
where $Z_i = \x{\lambda_r} \cdots \x{\lambda_{i+1}} \x{\lambda_{i-1}} \cdots \x{\lambda_1} \y{\mu_2} \cdots \y{\mu_s}$ and $R_i$ is a linear combination of terms 
$$ \x{\phi_{r-1}} \cdots \x{\phi_{1}} \y{\psi_1} \cdots \y{\psi_{s-1}}, \quad \phi_1, \dots, \phi_{r-1}, \psi_1, \dots, \psi_{s-1} \in \C.$$
If instead $k < i \leqslant r$, then $A_i$ is a linear combination of terms
$$ \x{\lambda_r} \cdots \x{\lambda_{i+1}} \x\gamma \x{\lambda_{i-1}} \cdots \x{\lambda_1} \y{\mu_2} \cdots \y{\mu_s},$$
where $\RootFn(\gamma) = \RootFn(\lambda_i) - \RootFn(\mu_1) \in \PositiveRoots$, since $\RootFn(\mu_1) \leqslant \RootFn(\lambda_1) \leqslant \RootFn(\lambda_i).$

Note that $ \Degree{\h}{\Proj(Z)} \leqslant \text{MAX}\set{\Degree{\h}{\Proj (A_i)}}$.
Consider now each of the three parts of the claim.

\textbf{Part \ref{ShapLemmaOne}.}  For $1 \leqslant i \leqslant k$,
\begin{equation}\label{EqualCase}
\Proj (A_i) = \Proj (Z_i) \LieBrac{\x{\lambda_i}}{\y{\mu_1}} + \Proj (R_i),
\end{equation}
since $\Proj|_{\UEA{\g}^0}$ is an algebra homomorphism.
By part \ref{ShapLemmaOne} of the inductive hypothesis,
$$ \Degree{\h}{\Proj (Z_i)}, \ \Degree{\h}{\Proj (R_i)} \leqslant r-1,\  s-1,$$
and so $\Degree{\h}{\Proj (A_i)} \leqslant r,s$.
For $k < i \leqslant r$, again by part \ref{ShapLemmaOne} of the inductive hypothesis, $\Degree{\h}{\Proj (A_i)} \leqslant r,\ s-1$.
Hence $\Degree{\h}{\Proj (Z)} \leqslant r,s$, and so part \ref{ShapLemmaOne} holds.

\textbf{Part \ref{ShapLemmaTwo}.}
Suppose that $r=s$, and let $\alpha \in \PositiveRoots$ be such that $\Length{\lambda^\alpha} \ne \Length{\mu^\alpha}$.
For $1 \leqslant i \leqslant k$, consider $\Proj (A_i)$ by equation (\ref{EqualCase}).
By part \ref{ShapLemmaOne} of the inductive hypothesis, 
$$ \Degree{\h}{\Proj (R_i)} \leqslant r-1 < r,$$
and so it remains only to consider $\Proj (Z_i)$.
Write $\lambda'$ (respectively, $\mu'$) for the partition consisting of the components of $\lambda$ (respectively, $\mu$) except for $\lambda_i$ (respectively, $\mu_1$).
Then $\Length{\lambda'} = \Length{\mu'}$ and $\Length{\lambda'^\alpha} \ne \Length{\mu'^\alpha}$.
Therefore, by part \ref{ShapLemmaTwo} of the inductive hypothesis, $\Degree{\h}{\Proj (Z_i)} < r-1 $; hence $\Degree{\h}{\Proj (A_i)} < r$.
For $k < i \leqslant r$, part \ref{ShapLemmaOne} of the inductive hypothesis implies that $\Degree{\h}{\Proj (A_i)} \leqslant s-1 < r$.
Therefore $\Degree{\h}{\Proj (Z)} < r$, as required.

\textbf{Part \ref{ShapLemmaThree}.}
Suppose that $r=s$ and that $\Length{\lambda^\alpha} = \Length{\mu^\alpha}$ for all $\alpha \in \PositiveRoots$.
Observe that for $k < i \leqslant r$, part \ref{ShapLemmaOne} of the inductive hypothesis implies that 
$$\Degree{\h}{\Proj (A_i)} \leqslant s-1 < r;$$
and that for $1 \leqslant i \leqslant k$, by the same, 
$$\Degree{\h}{\Proj (R_i)} \leqslant r-1 < r.$$
Therefore, the terms $\Proj (A_i)$ for $k < i \leqslant r$ and the terms $\Proj (R_i)$ for $1 \le i \le k$ can not contribute to the degree-$r$ component of $\Proj (Z)$; thus the degree-$r$ component of $\Proj (Z)$ is the degree-$r$ component of 
$$\sum_{i=1}^k \Proj(Z_i)\LieBrac{\x{\lambda_i}}{\y{\mu_1}}.$$
As $Z_i$ satisfies the conditions of part \ref{ShapLemmaThree} of the inductive hypothesis, for $1 \leqslant i \leqslant k$, the formula follows.
\end{proof}

\end{subsection}
\end{section}

\begin{section}{Truncated Current Lie Algebras}\label{TCLASection}
Let $(\g, \hZero, \h, \g_{+}, \Invol)$ be a Lie algebra with triangular decomposition, and let $\C$ denote a set parameterizing a root-basis for $\g_{+}$.
Fix a positive integer $\Nilp$, and let 
$$\gHat = \g \otimes \K[\Indet{t}] / \Indet{t}^{\Nilp+1} \K[\Indet{t}]$$
denote the associated truncated current Lie algebra with the triangular decomposition of Example \ref{TrunCurrentExample}.
The integer $\Nilp$ is the \NewTerm{nilpotency index} of $\gHat$.
Let $\CHat = \C \times \set{0, \dots, \Nilp}$.
Then $\CHat$ parameterizes a basis for $\gHat_{+}$ consisting of $\hZero$-weight vectors of homogeneous degree in $\Indet{t}$, via
$$ \CHat \ni \quad \gamma \quad \leftrightarrow \quad \x{\gamma} \quad \in \gHat_{+},$$
where $\x{\gamma} = \x{\tau} \otimes \Indet{t}^d$ if $\gamma = (\tau, d) \in \CHat$.
Define 
$$ \RootFn : \CHat \to \PositiveRoots, \qquad \DegreeFn{\Indet{t}}  : \CHat \to \set{0, \dots, \Nilp}$$ 
via $\x\gamma \in \RootSpace{\g}{\RootFn(\gamma)} \otimes \Indet{t}^{\DegreeFn{\Indet{t}}(\gamma)}$ for all $\gamma \in \CHat$.
Order the basis $\set { \x\gamma | \gamma \in \CHat}$ of $\gHat_{+}$ by fixing an arbitrary linearization of the partial order by increasing homogeneous degree in $\Indet{t}$, i.e. so that
$$ \DegreeFn{\Indet{t}} (\gamma) < \DegreeFn{\Indet{t}} (\gamma') \quad \Implies \quad \x\gamma < \x{\gamma'}, \qquad \gamma, \gamma' \in \CHat.$$
As per Subsection \ref{PartitionsSubsection}, the PBW basis monomials of $\UEA{\gHat_{+}}$ with respect to this ordered basis are parameterized by a collection $\Ptns$ of partitions.
Partitions here are (finite) multisets with elements from $\CHat$.
For any $\chi \in \QPlus$, let
$$ \Ptns_\chi = \set { \lambda \in \Ptns | \RootFn (\lambda) = \chi }.$$
For any $0 \leqslant d \leqslant \Nilp$ and $\lambda \in \Ptns$, define
$$ \lambda^d = \set{ \gamma \in \lambda | \Degree{\Indet{t}}{\gamma} = d};$$
$\lambda$ is homogeneous of degree-$d$ in $\Indet{t}$ if $\lambda = \lambda^d$.
The ordering of the basis $\gHat_{+}$ is such that for all $\lambda \in \Ptns$, 
$$ \x\lambda = \x{\lambda^0} \x{\lambda^1} \cdots \x{\lambda^\Nilp}, \qquad \y\lambda = \y{\lambda^\Nilp} \cdots \y{\lambda^1} \y{\lambda^0}.$$
For any $\Lambda \in \LinearDual\hHat$ and $0 \leqslant d \leqslant \Nilp$, let $\Lambda_d \in \LinearDual\h$ be given by
$$ \Functional{\Lambda_d}{h} = \Functional\Lambda{ h \otimes \Indet{t}^d}, \qquad h \in \h.$$
\begin{subsection}{The Shapovalov form}
As in Section \ref{HWTSection}, there is a decomposition
$$ \UEA{\gHat} = \UEA{\hHat} \oplus \{ \gHat_{-} \UEA{\gHat} + \UEA{\gHat} \gHat_{+} \},$$
as a direct sum of two-sided $\UEA{\gHat}^0$-modules.
Denote by $\Proj : \UEA{\gHat} \to \UEA{\hHat}$ the projection onto the first summand, parallel to the second.
Let $\ShapForm: \UEA{\gHat} \times \UEA{\gHat} \to \UEA{\hHat}$ denote the Shapovalov form, and write  $\ShapForm_\chi$ for the restriction of $\ShapForm$ to the subspace $\WeightSpace{\UEA{\gHat_{-}}}{-\chi}$, $\chi \in \QPlus$.

The algebra $\UEA\gHat$ is graded by total degree in the indeterminate $\Indet{t}$,
$$ \UEA\gHat = \bigoplus_{m \geqslant 0} \GradedComp{\UEA\gHat}{m}, \qquad 
\GradedComp{\UEA\gHat}{m} = \Span \set{ (x_1 \otimes \Indet{t}^{d_1}) \cdots (x_k \otimes \Indet{t}^{d_k}) | \sum_{i=1}^k d_i = m, \ \ k \geqslant 0}.$$
For any subspace $V \subset \UEA\gHat$, let
$$ \GradedComp{V}{m} = \GradedComp{\UEA\gHat}{m} \cap V, \qquad m \geqslant 0,$$
and call $V$ \NewTerm{graded} in $\Indet{t}$ if $V = \bigoplus_{m \geqslant 0} V_m$.
The subalgebras $\UEA{\gHat_{+}}$, $\UEA{\gHat_{-}}$, and $\UEA{\hHat}$ are graded in $\Indet{t}$.
\begin{lemma}\label{ProjPreservesTotalDegreeLemma}
For any $m \geqslant 0$, 
$ \Proj ( \GradedComp{\UEA\gHat}{m}) \subset \GradedComp{\UEA\hHat}{m}$.
\end{lemma}
\begin{proof}
The spaces $\gHat_{-} \UEA\gHat$ and $\UEA\gHat \gHat_{+}$ are graded in $\Indet{t}$; hence so is the sum 
$$
\{ \gHat_{-} \UEA{\gHat} + \UEA{\gHat} \gHat_{+} \}.
$$
Therefore,
$$ 
\GradedComp{\UEA{\gHat}}{m} = \GradedComp{\UEA{\hHat}}{m} \oplus \GradedComp{\{ \gHat_{-} \UEA{\gHat} + \UEA{\gHat} \gHat_{+} \}}{m},
$$
for any $m \geqslant 0$.
\end{proof}

\begin{example}\label{SL3PartitionsExample}
Let $\g = \SpecialLinear{3}$, and recall the notation of Example \ref{SL3Example}.
Let $\Nilp=1$, so that $\gHat = \g \oplus ({\g \otimes \Indet{t}})$.
Write $\C = \PositiveRoots$ and $\CHat = \C \times \set{0,1}$.
Then $\gHat_{+}$ has a basis parameterized by $\CHat$:
$$\CHat \ni \quad (\alpha, d) \quad \leftrightarrow \quad \x\alpha \otimes \Indet{t}^d \quad \in \gHat_{+}.$$
Let $\chi = \upalpha_1 + \upalpha_2$. Then $\Ptns_\chi$ consists of the six partitions
\begin{equation}\label{SL3Partitions}
\begin{array}{@{\extracolsep{1cm}}ccc}
\set{ (\upalpha_1, 0), (\upalpha_2, 0)},  &\set{ (\upalpha_1 + \upalpha_2, 0) },  &\set{ (\upalpha_1, 0), (\upalpha_2, 1) },\\
\set{ (\upalpha_1, 1), (\upalpha_2, 0)},  &\set{ (\upalpha_1 + \upalpha_2, 1) },  &\set{ (\upalpha_1, 1), (\upalpha_2, 1) }.
\end{array}
\end{equation}
Order the set $\set{\x{\gamma} | \gamma \in \CHat}$ by the enumeration
$$ \x{\upalpha_1}\otimes \Indet{t}^0, \quad \x{\upalpha_1 + \upalpha_2}\otimes \Indet{t}^0, \quad \x{\upalpha_2}\otimes \Indet{t}^0, \quad \x{\upalpha_1}\otimes \Indet{t}^1, \quad \x{\upalpha_1 + \upalpha_2}\otimes \Indet{t}^1, \quad \x{\upalpha_2}\otimes \Indet{t}^1.$$
Then the PBW basis monomials of $\WeightSpace{\UEA{\gHat_{-}}}{-\chi}$ corresponding to the partitions (\ref{SL3Partitions}) are, respectively,
\begin{equation}\label{sl3basis}
\begin{array}{@{\extracolsep{1cm}}ccc}
{\y{\upalpha_2} \otimes \Indet{t}^0} \cdot {\y{\upalpha_1} \otimes \Indet{t}^0}, & {\y{\upalpha_1 + \upalpha_2} \otimes \Indet{t}^0}, & {\y{\upalpha_2} \otimes \Indet{t}^1} \cdot {\y{\upalpha_1} \otimes \Indet{t}^0},\\
{\y{\upalpha_1} \otimes \Indet{t}^1} \cdot {\y{\upalpha_2} \otimes \Indet{t}^0}, & {\y{\upalpha_1 + \upalpha_2} \otimes \Indet{t}^1}, & {\y{\upalpha_2} \otimes \Indet{t}^1} \cdot {\y{\upalpha_1} \otimes \Indet{t}^1}.
\end{array}
\end{equation}
For notational convenience, write $\HT{\upalpha_i}{j} = \slh{\upalpha_i} \otimes \Indet{t}^j$, for $i=1,2$ and $j=0,1$. 
The restriction $\ShapForm_\chi$ of the Shapovalov form, expressed as a matrix with respect to the ordered basis (\ref{sl3basis}), appears below.  
$$
\begin{pmatrix}
\HT{\upalpha_1 + \upalpha_2}{0} + \HT{\upalpha_1}{0}	& \HT{\upalpha_1}{0} 	& \HT{\upalpha_1}{0} \HT{\upalpha_2}{1} + \HT{\upalpha_1}{1}	& \HT{\upalpha_1}{1} (\HT{\upalpha_2}{0} + 2) 	& \HT{\upalpha_1}{1} 	& \HT{\upalpha_1}{1} \HT{\upalpha_2}{1}	\\

\HT{\upalpha_1}{0}		& \HT{\upalpha_1 + \upalpha_2}{0} 	& \HT{\upalpha_1}{1}	& - \HT{\upalpha_2}{1}	& \HT{\upalpha_1 + \upalpha_2}{1}  & 0 \\

\HT{\upalpha_1}{0} \HT{\upalpha_2}{1} + \HT{\upalpha_1}{1}		& \HT{\upalpha_1}{1}	& 0 	& \HT{\upalpha_1}{1} \HT{\upalpha_2}{1} 	& 0	& 0 \\

\HT{\upalpha_1}{1} (\HT{\upalpha_2}{0} + 2) 	& - \HT{\upalpha_2}{1} 	& \HT{\upalpha_1}{1} \HT{\upalpha_2}{1}	& 0	& 0 	& 0 \\

\HT{\upalpha_1}{1}	& \HT{\upalpha_1 + \upalpha_2}{1} 	& 0 & 0 & 0 & 0 \\

\HT{\upalpha_1}{1} \HT{\upalpha_2}{1} & 0 & 0 & 0 & 0 & 0
\end{pmatrix}
$$
This is an elementary calculation using the commutation relations.
Observe that this matrix is triangular, and that in particular the determinant (the Shapovalov determinant at $\chi$) must be the product of the diagonal entries, viz.,
\begin{equation}\label{sl3determinant}
(\slh{\upalpha_1} \otimes \Indet{t}^1)^4 (\slh{\upalpha_2} \otimes \Indet{t}^1)^4 (\slh{\upalpha_1 + \upalpha_2} \otimes \Indet{t}^1)^2,
\end{equation}
up to sign.
This provides a criterion for the existence of primitive vectors in the weight space $\Lambda|_{\hZero} - \chi$ of a Verma module $\Verma{\Lambda}$, $\Lambda \in \LinearDual{\hHat}$.
We shall prove that the Shapovalov determinant always lies in $\SymAlg{\h \otimes \Indet{t}^\Nilp}$, and that for $\g$ a semisimple finite-dimensional Lie algebra, the factors of the Shapovalov determinant are the analogues of those of (\ref{sl3determinant}).
\end{example}

\begin{example}\label{VirasoroPartitionsExample1}
Let $\g$ be the Virasoro/Witt algebra, and adopt the notation of Example \ref{VirasoroExample}.
Let $\Nilp=1$, and $\chi = 2 \VirRoot$.
Write $\C = \PositiveRoots$ and $\CHat = \C \times \set{0,1}$.
Then $\Ptns_\chi$ consists of the five partitions
\begin{equation}\label{VirasoroPartitions1}
\begin{array}{ccccc}
\set{ (\VirRoot, 0), (\VirRoot, 0)},  &\set{ (2\VirRoot, 0) },  &\set{ (\VirRoot, 0), (\VirRoot, 1) }, &\set{ (2\VirRoot, 1) }, & \set{ (\VirRoot, 1), (\VirRoot, 1) }.
\end{array}
\end{equation}
Order the basis
$$ \set{\x{\gamma} | \gamma \in \CHat} = \set{\VirL_m \otimes \Indet{t}^d | m>0, \ \ d=0,1}$$
for $\gHat_{+}$ firstly by increasing degree $d$, and secondly by increasing index $m$.
Then the PBW basis monomials of $\WeightSpace{\UEA{\gHat_{-}}}{-\chi}$ corresponding to (\ref{VirasoroPartitions1}) are, respectively,
\begin{equation}\label{VirasoroBasisOne}
\begin{array}{ccccc}
({\VirL_{-1} \otimes \Indet{t}^0})^2,  &{\VirL_{-2} \otimes \Indet{t}^0}, &{\VirL_{-1} \otimes \Indet{t}^1} \cdot {\VirL_{-1} \otimes \Indet{t}^0}, &{\VirL_{-2} \otimes \Indet{t}^1}, &  ({\VirL_{-1} \otimes \Indet{t}^1})^2 .
\end{array}
\end{equation}
Write $\Vir{m}{i} = (2m \VirL_0 + \psi(m) \VirC) \otimes \Indet{t}^ i$, for $m,i \geqslant 0$. 
The matrix of $\ShapForm_\chi$, expressed with respect to the ordered basis (\ref{VirasoroBasisOne}), appears below.
$$
\begin{pmatrix}
2 \Vir{1}{0} ( \Vir{1}{0} + 1)	& 3 \Vir{1}{0}	& 2 \Vir{1}{1} (\Vir{1}{0} + 1)	& 3 \Vir{1}{1}	& 2 (\Vir{1}{1})^2 \\

3 \Vir{1}{0}	& \Vir{2}{0}	& 3 \Vir{1}{1}	& \Vir{2}{1}	& 0 \\

2 \Vir{1}{1} (\Vir{1}{0} + 1)	& 3 \Vir{1}{1}	& \Vir{1}{1}^2	& 0 & 0 \\

3 \Vir{1}{1}	& \Vir{2}{1}	& 0	& 0	& 0 \\

2 (\Vir{1}{1})^2	& 0 & 0 & 0 & 0

\end{pmatrix}
$$
Hence the Shapovalov determinant at $\chi$ is given by
$ \Det \ShapForm_\chi = 4 \Vir{1}{1}^6 \Vir{2}{1}^2$.
\end{example}

\begin{example}\label{VirasoroPartitionsExample2}
Let $\g$ be the Virasoro/Witt algebra, and adopt the notation of Example \ref{VirasoroExample}.
Let $\Nilp=2$, and $\chi = 2 \VirRoot$.
Write $\C = \PositiveRoots$ and $\CHat = \C \times \set{0,1,2}$.
Then $\Ptns_\chi$ consists of the nine partitions
\begin{equation}\label{VirasoroPartitions2}
\begin{array}{@{\extracolsep{1cm}}lll}
\set{ (\VirRoot, 0), (\VirRoot, 0)},  &\set{ (2\VirRoot, 0) },  &\set{ (\VirRoot, 0), (\VirRoot, 1) },\\
\set{ (\VirRoot, 0), (\VirRoot, 2) }, &\set{ (2\VirRoot, 1) }, &\set{ (\VirRoot, 1), (\VirRoot, 1)}, \\
\set{ (\VirRoot, 1), (\VirRoot, 2)}, &\set{ (2\VirRoot, 2) }, &\set{ (\VirRoot, 2), (\VirRoot, 2) }. 
\end{array}
\end{equation}
Order the basis $\set{\x{\gamma} | \gamma \in \CHat}$ as per Example \ref{VirasoroPartitionsExample1}.
Then the PBW basis monomials of $\WeightSpace{\UEA{\gHat_{-}}}{-\chi}$ corresponding to (\ref{VirasoroPartitions2}) are, respectively,
\begin{equation}\label{VirasoroBasisTwo}
\begin{array}{@{\extracolsep{1cm}}lll}
({\VirL_{-1} \otimes \Indet{t}^0})^2,  &{\VirL_{-2} \otimes \Indet{t}^0},	 &{\VirL_{-1} \otimes \Indet{t}^1} \cdot {\VirL_{-1} \otimes \Indet{t}^0},\\
{\VirL_{-1} \otimes \Indet{t}^2} \cdot {\VirL_{-1} \otimes \Indet{t}^0}, & {\VirL_{-2} \otimes \Indet{t}^1}, &({\VirL_{-1} \otimes \Indet{t}^1})^2,\\
{\VirL_{-1} \otimes \Indet{t}^2} \cdot {\VirL_{-1} \otimes \Indet{t}^1}, & {\VirL_{-2} \otimes \Indet{t}^2},	  &({\VirL_{-1} \otimes \Indet{t}^2})^2 . \\
\end{array}
\end{equation}
The matrix of $\ShapForm_\chi$ with respect to the ordered basis (\ref{VirasoroBasisTwo}) appears on page \pageref{BigVirasoroMatrix}.
Notice that the matrix has six non-zero entries on the diagonal.
Hence there is no reordering of the basis (\ref{VirasoroBasisTwo}) that will render the matrix triangular.
\end{example}
\end{subsection}
\begin{subsection}{A Modification of the Shapovalov Form}\label{ModifiedFormSubsection}
As observed in Example \ref{VirasoroPartitionsExample2}, it is not always the case that the matrix for $\ShapForm_\chi$, $\chi \in \QPlus$, can be made triangular by an ordering of the chosen PBW monomial basis for $\WeightSpace{\UEA{\gHat_{-}}}{-\chi}$.
A further permutation of columns is necessary; this is performed by an involution $\Star{}$ on the partitions, and encapsulated in a modification $\ModShap$ of the Shapovalov form $\ShapForm$.
For any $\gamma = (\tau, d) \in \CHat$, write $\Star{\gamma} = (\tau, \Nilp-d) \in \CHat$, and for any $\lambda \in \Ptns$, write
$$ \Star{\lambda} = \set {\Star{\gamma} | \gamma \in \lambda}.$$
So $\Star{(\lambda^d)} = (\Star{\lambda})^{\Nilp-d}$ for all $\lambda \in \Ptns$ and all degrees $d$.
For any $\chi \in \QPlus$, let
$$ \ModShap_\chi : \WeightSpace{\UEA{\gHat_{-}}}{-\chi} \times \WeightSpace{\UEA{\gHat_{-}}}{-\chi} \to \UEA{\hHat} $$
be the bilinear form defined by
$$ \ModShap_\chi (\y{\lambda}, \y{\mu}) = \ShapForm_\chi (\y{\lambda}, \y{\Star{\mu}}), \qquad \lambda, \mu \in \Ptns_\chi.$$
Relative to any linear order of the basis $\set{\y{\lambda} | \lambda \in \Ptns_\chi}$ of $\WeightSpace{\UEA{\gHat_{-}}}{-\chi}$, the matrices of $\ModShap_\chi$ and $\ShapForm_\chi$ are equal after a reordering of columns determined by the involution $\Star{}$.
In particular, the determinants $\Det \ModShap_\chi$ and $\Det \ShapForm_\chi$ are equal up to sign.
\end{subsection}
\end{section}

\begin{landscape}
\vspace*{2cm}
$$
\left(
\begin{array}{ccccccccc}\label{BigVirasoroMatrix}
2 \Vir{1}{0} (\Vir{1}{0} + 1) & 3 \Vir{1}{0} & 2 \Vir{1}{1} (\Vir{1}{0} + 1) & 2 \Vir{1}{2}(\Vir{1}{0} + 1) & 3 \Vir{1}{1} & 2 (\Vir{1}{1}^2 + \Vir{1}{2}) & 2 \Vir{1}{1} \Vir{1}{2} & 3 \Vir{1}{2} & 2 \Vir{1}{2}^2  \\

3 \Vir{1}{0} & \Vir{2}{0} & 3 \Vir{1}{1} & 3 \Vir{1}{2} & \Vir{2}{1} & 3 \Vir{1}{2} & 0 & \Vir{2}{2} & 0 \\

2 \Vir{1}{1} (\Vir{1}{0} + 1) & 3 \Vir{1}{1} & \Vir{1}{2} (\Vir{1}{0} + 2) + \Vir{1}{1}^2 & \Vir{1}{1} \Vir{1}{2} & 3 \Vir{1}{2} & 2 \Vir{1}{1} \Vir{1}{2} & \Vir{1}{2}^2 & 0 & 0 \\

2 \Vir{1}{2} (\Vir{1}{0} + 1) & 3 \Vir{1}{2} & \Vir{1}{1} \Vir{1}{2} & \Vir{1}{2}^2 & 0 & 0 & 0 & 0 & 0 \\

3 \Vir{1}{1} & \Vir{2}{1} & 3 \Vir{1}{2} & 0 & \Vir{2}{2} & 0 & 0 & 0 & 0 \\

2 (\Vir{1}{1}^2 + \Vir{1}{2}) & 3 \Vir{1}{2} & 2 \Vir{1}{1} \Vir{1}{2} & 0 & 0 & 2 \Vir{1}{2}^2 & 0 & 0 & 0 \\

2 \Vir{1}{1} \Vir{1}{2} & 0 & \Vir{1}{2}^2 & 0 & 0 & 0 & 0 & 0 & 0 \\

3 \Vir{1}{2} & \Vir{2}{2} & 0 & 0 & 0 & 0 & 0 & 0 & 0 \\

2 \Vir{1}{2}^2 & 0 & 0 & 0 & 0 & 0 & 0 & 0 & 0 
\end{array}
\right)
$$
\begin{center}
Matrix of the Shapovalov form $\ShapForm_\chi$ for the Virasoro/Witt truncated current Lie algebra. 
$$ \chi = 2 \VirRoot, \quad \Nilp =2, \quad \Vir{m}{i} := (2m \VirL_{0} + \psi (m) \VirC) \otimes \Indet{t}^i, \quad m,i \geqslant 0.$$
\end{center}
\end{landscape}

\begin{section}{Decomposition of the Shapovalov Form}\label{BlockDecompositionSection}
Throughout this section, let $(\g, \hZero, \h, \g_{+}, \Invol)$ denote a Lie algebra with triangular decomposition, and let $\gHat$ denote the associated truncated current Lie algebra of nilpotency index $\Nilp$.
Let $\MatrixSpace$ denote the collection of all two-dimensional arrays of non-negative integers with rows indexed by $\PositiveRoots$ and columns indexed by $\set{0, \dots, \Nilp}$, with only a finite number of non-zero entries.
For any $\chi \in \QPlus$, let
$$ \MatrixSpace_\chi = \set { L \in \MatrixSpace | \chi = \sum_{\alpha \in \PositiveRoots} \sum_{0 \leqslant d \leqslant \Nilp} L_{\alpha, d}\  \alpha }.$$
The entries of an array in $\MatrixSpace_\chi$ specify the multiplicity of each positive root in each homogeneous degree component of a partition of $\chi$, \ie 
$$ \lambda \in \Ptns_\chi \quad \iff \quad ( \ \Length{\lambda^{\alpha, d}} \ )_{\genfrac{}{}{0pt}{}{\alpha \in \PositiveRoots,}{0 \leqslant d \leqslant \Nilp} } \in \MatrixSpace_\chi.$$
Let 
$$\Ptns_{L} = \set{\lambda \in \Ptns | \Length{\lambda^{\alpha, d}} = L_{\alpha, d}, \ \text{for all} \ \alpha \in \PositiveRoots, \ 0 \leqslant d \leqslant \Nilp}.$$
Then for any $L \in \MatrixSpace$, $\Ptns_L$ is a non-empty finite set; if the root spaces of $\g$ are one-dimensional, then $\Ptns_L$ is a singleton.
The set $\MatrixSpace_\chi$ parameterizes a disjoint union decomposition of the set $\Ptns_\chi$:
\begin{equation}\label{PtnDisjointUnion}
\Ptns_\chi = \bigsqcup_{L \in \MatrixSpace_\chi} \Ptns_{L}.
\end{equation}
For any $S \subset \Ptns$, let
$$ \Span(S) = \Span_\K \set{ \y\lambda | \lambda \in S },$$
so that, for example, $\Span(\Ptns) = \UEA{\gHat_{-}}$ and $\Span(\Ptns_\chi) = \WeightSpace{\UEA{\gHat_{-}}}{- \chi}$ for any $\chi \in \QPlus$.
For any $\chi \in \QPlus$, the decomposition (\ref{PtnDisjointUnion}) of $\Ptns_\chi$ defines a decomposition of $\WeightSpace{\UEA{\gHat_{-}}}{-\chi} = \Span(\Ptns_\chi)$:
\begin{equation*}
\WeightSpace{\UEA{\gHat_{-}}}{-\chi} = \bigoplus_{L \in \MatrixSpace_\chi} \Span(\Ptns_L).
\end{equation*}
We construct an ordering of the set $\MatrixSpace_\chi$ and show that, relative to this ordering, any matrix expression of the modified Shapovalov form $\ModShap_{\chi}$ for $\gHat$ is block-upper-triangular (\cf Theorem \ref{TriangulationTheorem}).
The following Corollary, immediate from Theorem \ref{TriangulationTheorem}, provides a multiplicative decomposition of the Shapovalov determinant, and is the most important result of this section.
\begin{corollary}\label{DetBlockCorollary}
Let $\chi \in \QPlus$.  Then
$$ \Det{\ModShap_\chi} = \prod_{L \in \MatrixSpace_\chi} \Det{\ModShap_\chi |_{\Span(\Ptns_L)}}.$$
\end{corollary}

\begin{subsection}{An order on $\MatrixSpace_{\chi}$}
Write $\N = \set{0,1,2, 3, \dots }$ for the set of natural numbers, and fix an arbitrary linearization of the partial order on $\QPlus$.
If $X$ is a set with a linear order, write $\Reverse{X}$ for the set $X$ with the reverse order, \ie $x \leqslant y$ in $\Reverse{X}$ if and only if $x \geqslant y$ in $X$.
For example, the order on $\Reverse{\N}$ is such that $0$ is maximal.
Suppose that $(X_i)_{i \geqslant 1}$ is a sequence of linearly ordered sets, and let
$$ X = X_1 \times X_2 \times X_3 \times \cdots$$
denote the ordered Cartesian product.
The set $X$ carries an order $<_{X}$ defined by declaring, for all tuples $(x_i), (y_i) \in X$, that $(x_i) <_{X} (y_i)$ if and only if there exists some $m \geqslant 1$ such that $x_i = y_i$ for all $1 \leqslant i < m$, and $x_m < y_m$.
This order on $X$ is linear, and is called the \NewTerm{lexicographic order} (or \NewTerm{dictionary order}).
Fix an arbitrary enumeration of the countable set $\PositiveRoots \times \set{0, 1, \dots, \Nilp}$.
Consider $\MatrixSpace$ as a subset of the ordered Cartesian product of copies of the set $\N$ indexed by this enumeration.
Write $\OrderedMatrixSpace$ for the set $\MatrixSpace$ with the associated lexicographic order.
For any $L \in \MatrixSpace$, write
\begin{eqnarray*}
\RootFn(L) &=& (\sum_{\alpha \in \PositiveRoots} L_{\alpha, 0} \alpha, \sum_{\alpha \in \PositiveRoots} L_{\alpha, 1} \alpha, \dots, \sum_{\alpha \in \PositiveRoots} L_{\alpha, \Nilp} \alpha ) \quad \in (\Reverse{\QPlus})^{\Nilp + 1}, \\
\Length{L} &=& (\sum_{\alpha \in \PositiveRoots} L_{\alpha, 0}, \sum_{\alpha \in \PositiveRoots} L_{\alpha, 1}, \dots, \sum_{\alpha \in \PositiveRoots} L_{\alpha, \Nilp} ) \quad \in \Reverse{\N} \times \N^{\Nilp}. 
\end{eqnarray*}
For any $\chi \in \QPlus$, define a map
$$ \theta_\chi : \MatrixSpace_{\chi} \to (\Reverse{\QPlus})^{\Nilp + 1} \times (\Reverse{\N} \times \N^{\Nilp}) \times \OrderedMatrixSpace$$
by
$$ \theta_\chi (L) = (\RootFn(L), \Length{L}, L), \qquad L \in \MatrixSpace_{\chi}.$$
Consider the sets
$(\Reverse{\QPlus})^{\Nilp + 1}$ and $\Reverse{\N} \times \N^{\Nilp}$
to both carry lexicographic orders.
Thus the Cartesian product
\begin{equation}\label{OrderingSpace}
(\Reverse{\QPlus})^{\Nilp + 1} \times (\Reverse{\N} \times \N^{\Nilp}) \times \OrderedMatrixSpace
\end{equation}
carries a lexicographic order, and this order is linear.
For any $\chi \in \QPlus$, we consider the set $\MatrixSpace_{\chi}$ to carry the linear order defined by the injective map $\theta_\chi$ and the linearly ordered set (\ref{OrderingSpace}).
\end{subsection}

\begin{subsection}{Decomposition of the Shapovalov Form}
\begin{lemma}\label{DegreeLemma}
For any partitions $\lambda, \mu \in \Ptns$, 
$$ \x\lambda \y\mu = \x{\lambda^0} \y{\mu^\Nilp} \x{\lambda^1} \y{\mu^{\Nilp - 1}} \cdots \x{\lambda^\Nilp} \y{\mu^0}.$$
\end{lemma}
\begin{proof}
By choice of order for the basis $\set{\x\gamma | \gamma \in \CHat}$ of $\gHat_{+}$, and since $\y\mu = \Invol(\x\mu)$, by definition,
$$ \x\lambda \y\mu = \x{\lambda^0} \cdots \x{\lambda^\Nilp} \y{\mu^\Nilp} \cdots \y{\mu^0}.$$
As $\LieBrac{\x{\lambda^i}}{\y{\mu^j}} = 0$ if $i+j > \Nilp$, the claim follows.
\end{proof}

\begin{propn}\label{FormDecompPropn}
Suppose that $\lambda, \mu \in \Ptns_\chi$, $\chi \in \QPlus$, and further that $\Weight (\lambda^d) = \Weight (\mu^d)$, for all $0 \leqslant d \leqslant k$, for some $0 \leqslant k \leqslant \Nilp$.  Then
$$ \Form{\lambda}{\mu} = \prod_{0 \leqslant d \leqslant k}{\Form{\lambda^d}{\mu^d}}  \hspace{0.15em} \cdot \hspace{0.15em}
	\Form {\lambda'} {\mu'},$$
where $\lambda' = \bigcup_{k < d \leqslant \Nilp}{\lambda^d}$ and $\mu' = \bigcup_{k < d \leqslant \Nilp}{\mu^d}$.
\end{propn}
\begin{proof}
Under the hypotheses of the claim,
\begin{eqnarray*}
\Form{\lambda}{\mu} &=& \Proj (\x\lambda \y{\Star{\mu}}) \\ 
	&=& \Proj (\x{\lambda^0} \y{(\Star\mu)^\Nilp} \x{\lambda^1} \y{(\Star\mu)^{\Nilp - 1}} \cdots \x{\lambda^\Nilp} \y{(\Star\mu)^0}) \\
	&&\Note{by Lemma \ref{DegreeLemma}} \\
	&=& \Proj ( \x{\lambda^0} \y{\Star{(\mu^0)}} \x{\lambda^1} \y{\Star{(\mu^1)}} \cdots \x{\lambda^\Nilp} \y{\Star{(\mu^\Nilp)}}) \\
	&=& \prod_{0 \leqslant d \leqslant k}{\Proj(\x{\lambda^d} \y{\Star{\mu^{d}}})} \cdot \Proj ( \x{\lambda^{k+1}} \y{\Star{(\mu^{k+1})}} \cdots \x{\lambda^\Nilp} \y{\Star{(\mu^\Nilp)}}) \\
	&&\Note{since $\Proj |_{\UEA{\gHat}^0}$ is an algebra homomorphism} \\
	&=&\prod_{0 \leqslant d \leqslant k}{\Proj(\x{\lambda^d} \y{\Star{\mu^{d}}})}\cdot \Form {\lambda'}{\mu'} \\
	&&\Note{by Lemma \ref{DegreeLemma}} \\
	&=& (\prod_{0 \leqslant d \leqslant k}{\Form{\lambda^d}{\mu^d}}) \cdot \Form {\lambda'}{\mu'}. \mqed
\end{eqnarray*} 
\end{proof}

\begin{lemma}\label{VanishLemma}
Suppose that $\lambda, \mu \in \Ptns$ are partitions of homogeneous degree $d$.
\begin{enumerate}
\item\label{VanishLemmaPart1} If $d = 0$ and $\Length{\lambda} < \Length{\mu}$, or if $d > 0$ and $\Length{\lambda} > \Length{\mu}$, then $\Form{\lambda}{\mu} = 0$.
\item\label{VanishLemmaPart2} If $\Length{\lambda} = \Length{\mu}$ and $\Length{\lambda^{\alpha}} \ne \Length{\mu^{\alpha}}$ for some $\alpha \in \PositiveRoots$, then $\Form{\lambda}{\mu} = 0$.
\end{enumerate}
\end{lemma}
\begin{proof}
This Lemma follows essentially from Lemma \ref{ShapLemma}, applied to the Lie algebra with triangular decomposition $(\gHat, \hZero, \hHat, \gHat_{+}, \Invol)$.
Let $\lambda, \mu \in \Ptns$ be partitions of homogeneous degree $d$.
Since 
$$ \x{\lambda} \y{\Star{\mu}} \in \UEA{\gHat}_{\Length{\lambda}d + \Length{\mu}(\Nilp - d)},$$
it follows from Lemma \ref{ProjPreservesTotalDegreeLemma} that
\begin{equation}\label{VanishEqn1}
\Form{\lambda}{\mu} \in \UEA{\hHat}_{\Length{\lambda}d + \Length{\mu}(\Nilp - d)}.
\end{equation}
On the other hand,
\begin{equation}\label{VanishEqn2}
\Degree{\hHat}{\Form{\lambda}{\mu}} \leqslant \Length{\lambda}, \ \Length{\mu}
\end{equation}
by Lemma \ref{ShapLemma}.
Therefore, if $\Form{\lambda}{\mu} \ne 0$, and 
$$ \Form{\lambda}{\mu} \in \UEA{\hHat}_{m},$$
it must be that 
\begin{equation}\label{VanishEqn3}
m \leqslant \Length{\lambda}\Nilp \quad \text{and} \quad m \leqslant \Length{\mu} \Nilp,
\end{equation}
since the degree of $h \in \hHat$ in $\Indet{t}$ is at most $\Nilp$.
Combining (\ref{VanishEqn1}) and (\ref{VanishEqn3}), it follows that if $\Form{\lambda}{\mu} \ne 0$, then
\begin{equation}\label{VanishEqn4}
\Length{\lambda} d + \Length{\mu}(\Nilp - d) \leqslant \Length{\lambda}\Nilp,
\end{equation}
and
\begin{equation}\label{VanishEqn5}
\Length{\lambda} d + \Length{\mu} (\Nilp - d) \leqslant \Length{\mu} \Nilp.
\end{equation}
If $d=0$, then inequality (\ref{VanishEqn4}) becomes $\Length{\mu} \leqslant \Length{\lambda}$.
Hence, if $d=0$, and $\Length{\lambda} < \Length{\mu}$, then $\Form{\lambda}{\mu} = 0$.
If $d > 0$, then inequality (\ref{VanishEqn5}) yields $\Length{\lambda} \leqslant \Length{\mu}$.
Hence, if $d > 0$ and $\Length{\lambda} > \Length{\mu}$, it must be that $\Form{\lambda}{\mu} = 0$.
This proves part \ref{VanishLemmaPart1}.

Suppose now that $\Length{\lambda} = \Length{\mu} = r$, and that $\Length{\lambda^{\alpha}} \ne \Length{\mu^{\alpha}}$ for some $\alpha \in \PositiveRoots$.
Then, by Lemma \ref{ShapLemma}, the inequality (\ref{VanishEqn2}) becomes strict.
Hence, if $\Form{\lambda}{\mu} \ne 0$, then the inequalities (\ref{VanishEqn4}) and (\ref{VanishEqn5}) are also strict.
These both yield 
$ r \Nilp < r \Nilp$,
which is absurd.
Hence it must be that $\Form{\lambda}{\mu} = 0$, and part \ref{VanishLemmaPart2} is proven.
\end{proof}

\begin{lemma}\label{WeightLemma}
Suppose that $\nu \in \Q$ and $\nu \not \in \QPlus$.  Then
$\WeightSpace{\UEA{\gHat}}{\nu} \subset \gHat_{-} \UEA{\gHat}$.
\end{lemma}
\begin{proof}
Because $\gHat = \gHat_{-} \bigoplus (\hHat \bigoplus \gHat_{+})$, we have
$ \UEA{\gHat} = \UEA{\gHat_{-}} \bigotimes \UEA{ \hHat \bigoplus \gHat_{+} }$
by the PBW Theorem.
The set of all weights of the $\hZero$-module $\UEA{ \hHat \bigoplus \gHat_{+}}$ is precisely $\QPlus$, and so, for any $\nu \in \Q$,
\begin{equation}\label{WeightBreakdown}
\textstyle \WeightSpace{\UEA{\gHat}}{\nu} = \sum_{\chi \in \QPlus} \WeightSpace{\UEA{\gHat_{-}}}{\nu - \chi} \bigotimes \WeightSpace{\UEA{\hHat \bigoplus \gHat_{+}}}{\chi}.
\end{equation}
Suppose that $\nu \in \Q$ and $\nu \not \in \QPlus$.  Then, in particular,  $\nu - \chi \ne 0$, for any $\chi \in \QPlus$, and so
$$ \WeightSpace{\UEA{\gHat_{-}}}{\nu - \chi} \subset \gHat_{-} \UEA{\gHat}.$$
Hence $\UEA{\gHat}^\nu \subset \gHat_{-} \UEA{\gHat}$ by equation (\ref{WeightBreakdown}).
\end{proof}

\begin{theorem}\label{TriangulationTheorem}
Suppose that $\chi \in \QPlus$, and that $L,M \in \MatrixSpace_{\chi}$.
If $L > M$, then
$$ \Form{\lambda}{\mu} = 0,$$
for all $\lambda \in \Ptns_L$ and $\mu \in \Ptns_M$.
\end{theorem}
\begin{proof}
Suppose that $L, M \in \MatrixSpace_\chi$ and that $L > M$.
Then one of the following hold:
\begin{itemize}
\item $\RootFn(L) > \RootFn(M)$; or
\item $\RootFn(L) = \RootFn(M)$ and $\Length{L} > \Length{M}$; or
\item $\RootFn(L) = \RootFn(M)$, $\Length{L} = \Length{M}$ and $L > M$ in $\OrderedMatrixSpace$.
\end{itemize}
Let $\lambda \in \Ptns_L$ and let $\mu \in \Ptns_M$.

Suppose that $\RootFn(L) > \RootFn(M)$.
Then there exists $ 0 \leqslant l \leqslant \Nilp$ such that $\RootFn(\lambda^d) = \RootFn(\mu^d)$ for all $0 \leqslant d < l$, and $\RootFn(\lambda^l) > \RootFn(\mu^l)$ in $\Reverse{\QPlus}$, \ie $\RootFn(\lambda^l) < \RootFn(\mu^l)$ in $\QPlus$.
If $l > 0$, then Proposition \ref{FormDecompPropn} with $k = l - 1$ gives that
\begin{equation}\label{DecompThmEqn1}
\Form{ \lambda }{ \mu } = \theta \cdot \Form{ \lambda' }{ \mu' }
\end{equation}
for some $\theta \in \SymAlg{\hHat}$, where $\lambda' = \bigcup_{l \leqslant d \leqslant \Nilp}{\lambda^d}$ and $\mu' = \bigcup_{l \leqslant d \leqslant \Nilp}{\mu^d}$.
In the remaining case where $l = 0$, equation (\ref{DecompThmEqn1}) holds with $\lambda = \lambda'$, $\mu = \mu'$ and $\theta = 1$.
By Lemma \ref{DegreeLemma}, 
\begin{equation}\label{DecompThmEqn2}
\Form{ \lambda' }{ \mu' } = \Proj ( \x{\lambda^l} \y{ (\Star{\mu})^{\Nilp - l}} \cdots \x{\lambda^\Nilp} \y{ (\Star{\mu})^0}).
\end{equation}
Since $\RootFn((\Star{\mu})^{\Nilp - l}) = \RootFn(\mu^l)$, the monomial $\x{\lambda^l} \y{ (\Star{\mu})^{\Nilp - l} }$ has weight $\nu = \RootFn(\lambda^l) - \RootFn(\mu^l)$.
Now $\nu \not \in \QPlus$, since $\RootFn(\lambda^l) < \RootFn(\mu^l)$ in $\QPlus$, and so
$$
\x{\lambda^l} \y{ (\Star{\mu})^{\Nilp - l} } \in \gHat_{-} \UEA{\gHat}
$$
by Lemma \ref{WeightLemma}. Therefore
$$ \Form{ \lambda' }{ \mu' } = 0,$$
by equation (\ref{DecompThmEqn2}) and the definition of the projection $\Proj$.
Hence $\Form {\lambda}{\mu} = 0$ by equation (\ref{DecompThmEqn1}).

Suppose instead that $\RootFn(L) = \RootFn(M)$.  Then by Proposition \ref{FormDecompPropn},
$$
\Form{\lambda}{\mu} = \prod_{0 \leqslant d \leqslant \Nilp} \Form{\lambda^d}{\mu^d}.
$$
Suppose that $\Length{L} > \Length{M}$.
Then either $\Length{\lambda^d} < \Length{\mu^d}$, with $d= 0$, or $\Length{\lambda^d} > \Length{\mu^d}$ for some $0 < d \leqslant \Nilp$.
In either case, 
$$\Form{\lambda^d}{\mu^d} = 0,$$
by Lemma \ref{VanishLemma} part \ref{VanishLemmaPart1}, applied to the partitions $\lambda^d, \mu^d$.
Suppose that $\Length{L} = \Length{M}$ and that $L > M$ in $\OrderedMatrixSpace$.
Then 
$$ 
L_{\alpha, d} \ne M_{\alpha, d} \quad \text{for some} \quad  \alpha \in \PositiveRoots, \ 0 \leqslant d \leqslant \Nilp,
$$
so that $\Length{ (\lambda^d)^\alpha} \ne \Length{ (\mu^d)^\alpha}$.
Therefore, Lemma \ref{VanishLemma} part \ref{VanishLemmaPart2}, applied to the partitions $\lambda^d, \mu^d$, implies that $\Form{\lambda^d}{\mu^d} = 0$.
Hence $\Form{\lambda}{\mu} = 0$.
\end{proof}

\end{subsection}
\end{section}

\begin{section}{Values of the Shapovalov Form}\label{BlockSection}
Throughout this section, let $(\g, \hZero, \h, \g_{+}, \Invol)$ denote a Lie algebra with triangular decomposition, and let $\gHat$ denote the truncated current Lie algebra of nilpotency index $\Nilp$ associated to $\g$.
In Section \ref{BlockDecompositionSection}, the space $\WeightSpace{\UEA{\gHat_{-}}}{-\chi}$ is decomposed,
$$
\WeightSpace{\UEA{\gHat_{-}}}{-\chi} = \bigoplus_{L \in \MatrixSpace_\chi} \Span (\Ptns_L)
$$
and it is demonstrated that the determinant of the (modified) Shapovalov form $\ModShap_\chi$ on $\WeightSpace{\UEA{\gHat_{-}}}{-\chi}$ is the product of the determinants of the restrictions of $\ModShap_\chi$ to the spaces $\Span(\Ptns_L)$, $L \in \MatrixSpace_\chi$.
In this section, the restrictions $\ModShap|_{\Span(\Ptns_L)}$ are studied.
Firstly, the values of $\ModShap|_{\Span(\Ptns_L)}$ with respect to the basis $\y\lambda$, $\lambda \in \Ptns_L$, are calculated (\cf Proposition \ref{ExactValuePropn}).
This permits the recognition of $\ModShap|_{\Span(\Ptns_L)}$, in the case where $\g$ carries a non-degenerate pairing, as an $\SymAlg{\hHat}$-multiple of a non-degenerate bilinear form on $\Span(\Ptns_L)$ (\cf Theorem \ref{FormEqualityThm}).
The form on $\Span(\Ptns_L)$ is constructed as a symmetric tensor power of the non-degenerate form on $\g$.
\begin{subsection}{Values of the restrictions $\ModShap|_{\Span(\Ptns_L)}$}
Whenever $\lambda, \mu \in \Ptns$ and $\Length{\lambda} = \Length{\mu} = n$, let
$$ \SymPtnPair{\lambda}{\mu} = \sum_{\tau \in \Sym{n}} \prod_{1 \leqslant i \leqslant n} \LieBrac{ \x{\lambda_{\tau (i)}} }{ \y{\mu_i} } \quad \in \SymAlg{\hHat},$$
where $(\lambda_i)$ and $(\mu_i)$, $1 \leqslant i \leqslant n$ are arbitrary enumerations of $\lambda$ and $\mu$, respectively.
\begin{lemma}\label{SymPtnPairLemma}
Suppose that $\lambda, \mu \in \Ptns$ and $\Length{\lambda} = \Length{\mu}$.
\begin{enumerate}
\item\label{SymPtnPairLemma1} $\SymPtnPair{\lambda}{\mu} = \SymPtnPair{\mu}{\lambda};$
\item\label{SymPtnPairLemma2} if, in addition, $\lambda$ and $\mu$ are homogeneous of degree-$d$ in $\Indet{t}$, then
$$ \SymPtnPair{\Star{\lambda}}{\mu} = \SymPtnPair\lambda{\Star\mu}.$$
\end{enumerate}
\end{lemma}
\begin{proof}
Let $n = \Length{\lambda} = \Length{\mu}$, and choose some enumerations $(\lambda_i)$, $(\mu_i)$, $1 \leqslant i \leqslant n$ of $\lambda$ and $\mu$.
The anti-involution $\Invol$ point-wise fixes $\SymAlg{\hHat}$, and so $\Invol$ fixes $\SymPtnPair{\lambda}{\mu}$.
On the other hand,
\begin{eqnarray*}
\Invol (\SymPtnPair{\lambda}{\mu}) &=& \sum_{\tau \in \Sym{n}} \prod_{1 \leqslant i \leqslant n} \Invol (\LieBrac{ \x{\lambda_{\tau (i)}} }{ \y{\mu_i} }) \\
 &=& \sum_{\tau \in \Sym{n}} \prod_{1 \leqslant i \leqslant n}\LieBrac{ \x{\mu_i} }{ \y{\lambda_{\tau (i)}}} \\
 &=& \SymPtnPair{\mu}{\lambda}
\end{eqnarray*}
proving part \ref{SymPtnPairLemma1}.
Suppose that $\lambda$, $\mu$ are homogeneous of degree-$d$ in $\Indet{t}$. 
For each $1 \leqslant i \leqslant n$, let $\epsilon_i, \gamma_i \in \C$ be such that
$$ \lambda_i = (\epsilon_i, d), \qquad \mu_i = (\gamma_i, d).$$
Then
\begin{equation}\label{SymPtnStarEquation}
\SymPtnPair{\Star{\lambda}}{\mu} = \sum_{\tau \in \Sym{n}} \prod_{1 \leqslant i \leqslant n} \LieBrac
	{ \x{ \epsilon_{\tau (i)} } \otimes \Indet{t}^{\Nilp - d} }
	{ \y{ \gamma_i} \otimes \Indet{t}^d}.
\end{equation}
For any $1 \leqslant i \leqslant n$ and $\tau \in \Sym{n}$,
\begin{eqnarray*}
\LieBrac {\x{\epsilon_{\tau (i)}} \otimes \Indet{t}^{\Nilp - d}}{\y{\gamma_i} \otimes \Indet{t}^d} &=&  \LieBrac {\x{\epsilon_{\tau (i)}} }{\y{\gamma_i} } \otimes \Indet{t}^\Nilp \\
	&=& \LieBrac {\x{\epsilon_{\tau (i)}} \otimes \Indet{t}^{d}}{\y{\gamma_i} \otimes \Indet{t}^{\Nilp -d }},
\end{eqnarray*}
and hence $ \SymPtnPair{\Star{\lambda}}{\mu} = \SymPtnPair\lambda{\Star\mu}$ by equation (\ref{SymPtnStarEquation}), proving part \ref{SymPtnPairLemma2}.
\end{proof}
\begin{propn}\label{ExactValuePropn}
Suppose that $L \in \MatrixSpace$, and that $\lambda, \mu \in \Ptns_L$. 
Then
\begin{equation}\label{ExactValueFormula}
\Form{\lambda}{\mu} = \prod_{0 \leqslant d \leqslant \Nilp} \prod_{\alpha \in \PositiveRoots} \SymPtnPair{\lambda^{\alpha, d}}{\Star{(\mu^{\alpha, d})}}
\end{equation}
and $\Form{\lambda}{\mu} = \Form{\mu}{\lambda}$.
\end{propn}
\begin{proof}
Let $\lambda, \mu \in \Ptns_L$, $L \in \MatrixSpace$, and let $0 \leqslant d \leqslant \Nilp$.
Then 
$$\Length {\lambda^{\alpha, d}} = \Length {\mu^{\alpha, d}} = L_{\alpha,d}, \quad \alpha \in \PositiveRoots.$$
Write $l = \Length {\lambda^d} = \Length {\mu^d}$.
Then by Lemma \ref{ShapLemma}, applied to the Lie algebra with triangular decomposition $(\gHat, \hZero, \hHat, \gHat_{+}, \Invol)$,
$$ \Degree{\hHat}{\Form{\lambda^d}{\mu^d}} \leqslant l,$$
and the degree-$l$ component of $\Form{\lambda^d}{\mu^d}$ is given by 
\begin{equation}\label{ExactValuePropnEqn1}
\prod_{\alpha \in \PositiveRoots} \SymPtnPair{\lambda^{\alpha, d}}{ \Star{(\mu^{\alpha, d})}},
\end{equation}
since $\Form{\lambda^d}{\mu^d} = \Proj( \x{\lambda^d} \y{\mu^d {}^\star})$.
By Lemma \ref{ProjPreservesTotalDegreeLemma}, and since $ ld + l(\Nilp-d) = l\Nilp$, 
$$\Form{\lambda^d}{\mu^d} \in \GradedComp{\UEA\hHat}{l\Nilp}.$$ 
Therefore $\Degree{\hHat}{\Form{\lambda^d}{\mu^d}} \geqslant l$, since $\DegreeFn{\Indet{t}} \phi \leqslant \Nilp$ for any $\phi \in \hHat$; and so $\Form{\lambda^d}{\mu^d}$ is homogeneous of degree-$l$ in $\hHat$, and is equal to the expression (\ref{ExactValuePropnEqn1}).
By Proposition \ref{FormDecompPropn},
$$ \Form{\lambda}{\mu} = \prod_{0 \leqslant d \leqslant \Nilp}{\Form{\lambda^d}{\mu^d}},$$
and so the equation (\ref{ExactValueFormula}) follows.
The symmetry of $\ModShap|_{\Span(\Ptns_L)}$ follows from equation (\ref{ExactValueFormula}),
\begin{eqnarray*}
\Form{\lambda}{\mu} &=& \prod_{0 \leqslant d \leqslant \Nilp} \prod_{\alpha \in \PositiveRoots} \SymPtnPair{\lambda^{\alpha, d}}{\Star{(\mu^{\alpha, d})}} \\
	&=& \prod_{0 \leqslant d \leqslant \Nilp} \prod_{\alpha \in \PositiveRoots} \SymPtnPair{\Star{(\mu^{\alpha, d})}}{\lambda^{\alpha, d}} \\
	&=& \prod_{0 \leqslant d \leqslant \Nilp} \prod_{\alpha \in \PositiveRoots} \SymPtnPair{\mu^{\alpha, d}}{\Star{(\lambda^{\alpha, d})}} \\
	&=& \Form\mu\lambda
\end{eqnarray*}
and parts \ref{SymPtnPairLemma1} and \ref{SymPtnPairLemma2} of Lemma \ref{SymPtnPairLemma}.
\end{proof}
\end{subsection}

\begin{subsection}{Tensor powers of bilinear forms}
If $U$, $V$ are vector spaces, and $\phi : U \times V \to \K$ is a bilinear map, write
\begin{equation}\label{TensorUniversal}
\tilde{\phi} : \TensorProduct{U}{V} \to \K
\end{equation}
for the unique linear map such that $\tilde{\phi} ( u \otimes v) = \phi (u,v)$ for all $u \in U$, $v \in V$.
\begin{propn}\label{TensorFormPropn}
Suppose that $U$, $V$ are vector spaces with bilinear forms.
Then the vector space $U \bigotimes V$ carries a bilinear form defined by
\begin{equation}\label{TensorForm}
\VertForm{ u_1 \otimes v_1 } {u_2 \otimes v_2 } = \VertForm{u_1}{u_2} \VertForm{v_1}{v_2},
\end{equation}
for all $u_1, u_2 \in U$ and $v_1, v_2 \in V$.
Moreover, if the forms on $U$ and $V$ are non-degenerate, then so is the form on $U \bigotimes V$.
\end{propn}
\begin{proof}
Let $\phi : U \times U \to \K$, $\psi : V \times V \to \K$ denote the bilinear forms on $U$, $V$, respectively.
Let
$$
\nu : (\TensorProduct{U}{U}) \times ( \TensorProduct{V}{V} ) \to \K
$$
be given by
$$
\nu ( u_1 \otimes u_2, v_1 \otimes v_2) = \tilde{\phi}(u_1 \otimes u_2) \tilde{\psi}(v_1 \otimes v_2),
$$
for all $u_1, u_2 \in U$, $v_1, v_2 \in V$, where the maps $\tilde{\phi}, \tilde{\psi}$ are defined by (\ref{TensorUniversal}).
Then $\nu$ is bilinear, and so defines a linear map
$$ \tilde{\nu} : \TensorProduct{(\TensorProduct{U}{U}) }{ (\TensorProduct{V}{V})} \to \K $$
by (\ref{TensorUniversal}).  Since
$$ \TensorProduct{(\TensorProduct{U}{U}) }{ (\TensorProduct{V}{V})} \cong \TensorProduct{(\TensorProduct{U}{V}) }{(\TensorProduct{U}{V})},$$
the map $\tilde{\nu}$ may be considered as a bilinear form
$$
\VertForm{\cdot}{\cdot} : (\TensorProduct{U}{V}) \times (\TensorProduct{U}{V}) \to \K.
$$
Now if $u_1, u_2 \in U$, $v_1, v_2 \in V$, then
\begin{eqnarray*}
\VertForm{u_1 \otimes v_1}{u_2 \otimes v_2} &=& \tilde{\nu} (u_1 \otimes u_2 \otimes v_1 \otimes v_2) \\
	&=& \nu (u_1 \otimes u_2, v_1 \otimes v_2) \\
	&=& \tilde{\phi} (u_1 \otimes u_2) \tilde{\psi}(v_1 \otimes v_2) \\
	&=& \phi(u_1, u_2) \psi(v_1, v_2),
\end{eqnarray*}
and so this is the required bilinear form.
The non-degeneracy claim follows immediately from the definition (\ref{TensorForm}) of the form.
\end{proof}

For any vector space $U$ and non-negative integer $n$, write
$$ \textstyle \TensorPower{n}{U} = U \bigotimes \cdots \bigotimes U, \ENote{$n$ times}$$
for the space of homogeneous degree-$n$ tensors in $U$.
For any $u_i \in U$, $1 \leqslant i \leqslant n$, write 
$$ \textstyle \bigotimes_{i=1}^{n} u_i = u_1 \otimes \cdots \otimes u_n,$$
so that
$$ \TensorPower{n}{U} = \Span \set{ \textstyle \bigotimes_{i=1}^{n}{u_i}  | u_i \in U, \ 1 \leqslant i \leqslant n}.$$
Write
$$ u_1 \cdots u_n = \prod_{i=1}^{n} u_i = \frac{1}{\Factorial{n}} \sum_{\sigma \in \Sym{n}} {\textstyle \bigotimes_{i=1}^{n} u_{\sigma(i)}}
$$
for the symmetric tensor in $u_i \in U$, $1 \leqslant i \leqslant n$, and let
$$ \SymPower{n}{U} = \Span \set{ \prod_{i=1}^{n} u_i | u_i \in U, \ 1 \leqslant i \leqslant n}$$
denote the space of degree-$n$ symmetric tensors in $U$.
Let
$$ \TensorAlg{U} = \bigoplus_{n \ge 0} \TensorPower{n}{U}, \qquad \SymAlg{U} = \bigoplus_{n \ge 0} \SymPower{n}{U},$$
denote the tensor and symmetric algebras over $U$, respectively.
\begin{propn}\label{SymPowerPropn}
Suppose that $U$ is a vector space endowed with a bilinear form, and that $n \geqslant 0$.
Then $\SymPower{n}{U}$ carries a bilinear form defined by
$$ \VertForm{ \prod_{i=1}^{n} {u_i}}{\prod_{i=1}^{n} {v_i }} = \frac{1}{\Factorial{n}} \sum_{\tau \in \Sym{n}} \prod_{i=1}^{n} \VertForm{u_i}{v_{\tau (i)}},$$
for any $u_i, v_i \in U$, $1 \leqslant i \leqslant n$.
Moreover, if the form on $U$ is non-degenerate, then so is the form on $\SymPower{n}{U}$.
\end{propn}
\begin{proof}
Let $\text{A}(U)$ denote the two-sided ideal of $\TensorAlg{U}$ generated by the elements of the set
$$ \set{u_1 \otimes u_2 - u_2 \otimes u_1 | u_1, u_2 \in U}.$$
Then $\TensorAlg{U} = \SymAlg{U} \bigoplus \text{A}(U)$ is a direct sum of graded vector spaces.
Hence, for any $n \geqslant 0$, 
\begin{equation}\label{OrthogDecomp}
 \TensorPower{n}{U} = \SymPower{n}{U} \bigoplus \text{A}^{n}(U)
\end{equation}
is a direct sum of vector spaces, where $\text{A}^{n}(U)$ denotes the homogeneous degree-$n$ component of $\text{A}(U)$.
By Proposition \ref{TensorFormPropn}, the tensor power $\TensorPower{n}{U}$ carries a bilinear form defined by
\begin{equation}\label{TensorPowerForm}
\VertForm{ {\textstyle \bigotimes_{i=1}^{n} u_i} }{ {\textstyle \bigotimes_{i=1}^{n} v_i } } = \prod_{i=1}^{n} \VertForm{u_i}{ v_i}, \qquad u_i, v_i \in U, \quad 1 \leqslant i \leqslant n.
\end{equation}
Observe that for any $u_i, v_i \in U$, $1 \leqslant i \leqslant n$, 
$$ \sum_{\sigma \in \Sym{n}} \prod_{i=1}^{n} \VertForm{u_i}{ v_{\sigma (i)}}$$
is independent of the enumeration of the elements $v_1, \dots, v_n$.
It follows that the direct sum (\ref{OrthogDecomp}) is orthogonal, with respect to the bilinear form (\ref{TensorPowerForm}).
A form is defined on $\SymPower{n}{U}$ by restriction of the form on $\TensorPower{n}{U}$.
For any $u_i, v_i \in U$, $1 \leqslant i \leqslant n$, 
\begin{eqnarray*}
\VertForm{ \prod_{i=1}^{n} u_i }{ \prod_{i=1}^{n} v_i } &=& (\frac{1}{\Factorial{n}})^2 \sum_{\sigma \in \Sym{n}} \sum_{\tau \in \Sym{n}} \prod_{i=1}^{n} \VertForm{u_{\sigma (i)}}{ v_{\tau (i)}} \\
	&=& \frac{1}{\Factorial{n}} \sum_{\tau \in \Sym{n}} \prod_{i=1}^{n} \VertForm{u_i}{ v_{\tau(i)}}. 
\end{eqnarray*}
Hence $\SymPower{n}{U}$ carries the required bilinear form.
If the form on $U$ is non-degenerate, 
then by Proposition \ref{TensorFormPropn} the form on $\TensorPower{n}{U}$ is non-degenerate, and
since the sum (\ref{OrthogDecomp}) is orthogonal, the restriction of the form to $\SymPower{n}{U}$ is non-degenerate also.
\end{proof}
\end{subsection}

\begin{subsection}{Lie algebras with non-degenerate pairing}
A Lie algebra with triangular decomposition $(\g, \hZero, \h, \g_{+}, \Invol)$ is said to have \NewTerm{non-degenerate pairing} if for all $\alpha \in \Roots$, there exists a non-zero $\hElement{\alpha} \in \h$, and a non-degenerate bilinear form
$$ \PairForm{\cdot}{\cdot}{\alpha} : \RootSpace{\g}{\alpha} \times \RootSpace{\g}{\alpha} \to \K,$$
such that 
\begin{equation}\label{NonDegPairingEqn}
\LieBrac{x_1}{\Invol(x_2)} = \PairForm{x_1}{x_2}{\alpha} \hElement{\alpha},
\end{equation}
for all $x_1, x_2 \in \RootSpace{\g}{\alpha}$.

If $\g$ has a non-degenerate pairing, then for any $\alpha \in \Roots$, the space
$$ \LieBrac{\RootSpace\g\alpha}{\RootSpace\g{-\alpha}} = \LieBrac{\RootSpace\g{-\alpha}}{\RootSpace\g{\alpha}}$$
is one-dimensional, and so the elements $\hElement\alpha$ and $\hElement{-\alpha}$ can differ only by a non-zero scalar.
\begin{example}\label{KMLAPairingExample}
Let $\g$ be a symmetrizable Kac-Moody Lie algebra over $\K$ (\cf Example \ref{KMLAExample}), and let $\VertForm{\cdot}{\cdot}$ denote a standard bilinear form on $\g$ (as per \cite[page 20]{KacBook}).
The restriction of this form to $\h$ is non-degenerate.
Therefore, for any $\chi \in \LinearDual{\h}$, there exists a unique $\hElement{\chi} \in \h$ such that 
$$ \Functional{\chi}{h} = \VertForm{\hElement{\chi}}{h} \qquad h \in \h.$$
The map $\hElementMap : \LinearDual{\h} \to \h$ is a linear isomorphism.
For any $\alpha \in \Roots$, let
$$ \PairForm{\cdot}{\cdot}{\alpha} : \RootSpace{\g}{\alpha} \times \RootSpace{\g}{\alpha} \to \K,$$
be given by
$$ \PairForm{x_1}{x_2}{\alpha} = \VertForm{x_1}{\Invol (x_2)}, \qquad x_1, x_2 \in \RootSpace{\g}{\alpha}.$$
Then for any $\alpha \in \Roots$, the form $\PairForm{\cdot}{\cdot}{\alpha}$ is non-degenerate, and is such that equation (\ref{NonDegPairingEqn}) holds (see, for example, Theorem 2.2 of \cite{KacBook}).
Hence $\g$ carries a non-degenerate pairing.
\end{example}

\begin{example}\label{FinePairingExample}
Suppose that $\g$ is a Lie algebra with triangular decomposition, such that for any root $\alpha \in \Roots$,
$$ \Dim \RootSpace\g\alpha = \Dim \RootSpace\g{-\alpha} = 1,
\quad \text{and} \quad 
\LieBrac{\RootSpace\g\alpha}{\RootSpace\g{-\alpha}} \ne 0.$$
Then for each $\alpha \in \Roots$, we may choose an arbitrary non-zero
$$
\hElement\alpha \in \LieBrac{\RootSpace\g\alpha}{\RootSpace\g{-\alpha}}
$$
and let the form $\PairForm\cdot\cdot\alpha : \RootSpace\g\alpha \times \RootSpace\g\alpha \to \K$ be defined by equation (\ref{NonDegPairingEqn}).
\end{example}

\begin{example}\label{VirasoroPairingExample}
Let $\g$ denote the Virasoro algebra (\cf Example \ref{VirasoroExample}).
Let $\alpha \in \Roots$, and let $m$ be the non-zero integer such that $\alpha = m \VirRoot$.
Then
$$ \RootSpace\g\alpha = \K \VirL_m, \qquad \RootSpace\g{-\alpha} = \K \VirL_{-m},$$
and
$
\LieBrac{\VirL_m}{\VirL_{-m}} = 2m \VirL_0 + \psi(m) \VirC
$
is non-zero.
Therefore, by Example \ref{FinePairingExample}, $\g$ carries a non-degenerate pairing, with
$$ \hElement\alpha = 2m \VirL_0 + \psi(m) \VirC, \qquad \alpha = m \VirRoot.$$
\end{example}

\begin{example}\label{HeisenbergPairingExample}
The Heisenberg Lie algebra $\Heis$ carries a non-degenerate pairing (\cf Example \ref{HeisenbergExampleOne}).
Let $\alpha \in \Roots$, and let $m$ be the non-zero integer such that $\alpha = m \HeisRoot$.
Then
$$ \RootSpace\Heis\alpha = \K \HeisA_m, \qquad \RootSpace\Heis{-\alpha} = \K \HeisA_{-m},$$
and
$
\LieBrac{\HeisA_m}{\HeisA_{-m}} = m \HeisH 
$
is non-zero.
Therefore, by Example \ref{FinePairingExample}, $\Heis$ carries a non-degenerate pairing, with
$$ \hElement\alpha = m \HeisH, \qquad \alpha = m \HeisRoot.$$
\end{example}

Suppose that $(\g, \hZero, \h, \g_{+}, \Invol)$ is a Lie algebra with triangular decomposition and non-degenerate pairing, and let $\gHat$ denote the truncated current Lie algebra with nilpotency index $\Nilp$ associated to $\g$.
Non-degenerate bilinear forms are defined on the homogeneous degree components of the roots spaces of $\gHat$ in the following manner.
For all $\alpha \in \PositiveRoots$ and $0 \leqslant d \leqslant \Nilp$, define
a non-degenerate bilinear form $\PairForm{\cdot}{\cdot}{\alpha, d}$ on $\RootSpace{\g}{\alpha}\otimes \Indet{t}^d$ by
\begin{equation}\label{NondegFormExtn}
\PairForm{x_1 \otimes \Indet{t}^d}{x_2 \otimes \Indet{t}^d}{\alpha, d} = \PairForm{x_1}{x_2}{\alpha}, \qquad x_1, x_2 \in \RootSpace{\g}{\alpha}.
\end{equation}
For all $\alpha \in \PositiveRoots$ and $0 \leqslant d \leqslant \Nilp$, let 
$$ \CHat_{\alpha, d} = \set{ \gamma \in \CHat | \RootFn(\gamma) = \alpha, \quad \Degree{\Indet{t}}{\gamma} = d }.$$
\begin{lemma}\label{PairingLemma}
Let $\alpha \in \PositiveRoots$ and let $0 \leqslant d \leqslant \Nilp$.  Then,
$$ \LieBrac{ \x{\phi}}{\y{\Star{\psi}}} = \PairForm{ \x{\phi}}{\x{\psi}}{\alpha,d} \hElement{\alpha}\otimes\Indet{t}^\Nilp,$$
for all $\phi, \psi \in \CHat_{\alpha, d}$.
\end{lemma}
\begin{proof}
Let $\phi^\prime, \psi^\prime \in \C$ be such that $\phi = (\phi^\prime,d)$ and $\psi = (\psi^\prime,d)$.  Then
$$ \y{\Star{\psi}} = \y{\psi^\prime}\otimes \Indet{t}^{\Nilp - d} = \Invol ( \x{\psi^\prime}) \otimes \Indet{t}^{\Nilp - d},$$
and $\x{\phi} = \x{\phi^\prime} \otimes \Indet{t}^d$.  Therefore
\begin{eqnarray*}
\LieBrac{\x{\phi}}{\y{\Star{\psi}}} &=& \LieBrac{\x{\phi^\prime}}{\Invol(\x{\psi^\prime})}\otimes\Indet{t}^\Nilp \\
	&=& \PairForm{\x{\phi^\prime}}{\x{\psi^\prime}}{\alpha} \hElement{\alpha} \otimes \Indet{t}^\Nilp \\
	&=& \PairForm{\x{\phi}}{\x{\psi}}{\alpha, d} \hElement{\alpha} \otimes \Indet{t}^\Nilp,
\end{eqnarray*}
by equation (\ref{NondegFormExtn}).
\end{proof}
\end{subsection}
\begin{subsection}{Recognition of the restrictions $\ModShap|_{\Span(\Ptns_L)}$}
For any $L \in \MatrixSpace$, let
$$ \SymSpace{L} = \bigotimes_{0 \leqslant d \leqslant \Nilp} \bigotimes_{\alpha \in \PositiveRoots} \SymPower{L_{\alpha,d}}{\RootSpace{\g}{\alpha} \otimes \Indet{t}^d}.$$
The vector space $\SymSpace{L}$ has a basis parameterized by the partitions in $\Ptns_{L}$:
$$ \Ptns_{L} \ni \quad \lambda \quad \leftrightarrow \quad  \SymSpaceBasis{\lambda} \quad \in \SymSpace{L},$$
where, for all $\lambda \in \Ptns$,
$$ \SymSpaceBasis{\lambda} = \bigotimes_{0 \leqslant d \leqslant \Nilp} \bigotimes_{\alpha \in \PositiveRoots} \prod_{\gamma \in \lambda^{\alpha,d}} \x{\gamma}.$$

\begin{propn}\label{SymPtnPairValuePropn}
Let $\alpha \in \PositiveRoots$ and $0 \leqslant d \leqslant \Nilp$.  
If $\lambda, \mu \in \Ptns$ are partitions with components in $\CHat_{\alpha, d}$ such that $\Length{\lambda} = \Length{\mu} = k$, then
$$ \SymPtnPair{\lambda}{\Star{\mu}} = \Factorial{k} \hspace{0.2em} ( \hElement{\alpha} \otimes \Indet{t}^\Nilp)^k \hspace{0.2em} \VertForm{\SymSpaceBasis{\lambda}}{\SymSpaceBasis{\mu}}$$
where $\VertForm{\cdot}{\cdot}$ is the form on $\SymPower{k}{\RootSpace{\g}{\alpha}\otimes\Indet{t}^d}$ defined by the form on $\RootSpace{\g}{\alpha}\otimes\Indet{t}^d$ and Proposition \ref{SymPowerPropn}.
\end{propn}
\begin{proof}
The claim follows from Lemma \ref{PairingLemma} and the definition of the form on $\SymPower{k}{\RootSpace{\g}{\alpha}\otimes\Indet{t}^d}$.
Let $(\lambda_i)$ and $(\mu_i)$, $1 \leqslant i \leqslant k$ be any enumerations of $\lambda$ and $\mu$, respectively.  Then:
\begin{eqnarray*}
\SymPtnPair{\lambda}{\Star{\mu}} &=& \sum_{\tau \in \Sym{k}} \prod_{1 \leqslant i \leqslant k} \LieBrac{\x{\lambda_{\tau (i)}}}{\y{\Star{\mu_i}}} \\
	&=& \sum_{\tau \in \Sym{k}} \prod_{1 \leqslant i \leqslant k} \PairForm{\x{\lambda_{\tau(i)}}}{\x{\mu_i}}{\alpha, d} \hElement{\alpha} \otimes \Indet{t}^\Nilp \\
	&=& \Factorial{k} \hspace{0.2em} (\hElement{\alpha}\otimes\Indet{t}^\Nilp)^k \hspace{0.2em} \frac{1}{\Factorial{k}} \sum_{\tau \in \Sym{k}} \prod_{1 \leqslant i \leqslant k} \PairForm{\x{\lambda_{\tau(i)}}}{\x{\mu_i}}{\alpha, d} \\
	&=& \Factorial{k} \hspace{0.2em} (\hElement{\alpha}\otimes\Indet{t}^\Nilp)^k \hspace{0.2em} \VertForm{\SymSpaceBasis{\lambda}}{\SymSpaceBasis{\mu}}.\mqed
\end{eqnarray*}
\end{proof}

For any $L \in \MatrixSpace$, the vector spaces $\Span(\Ptns_{L})$ and $\SymSpace{L}$ are isomorphic by linear extension of the correspondence
$$ \y{\lambda} \quad \leftrightarrow \quad \SymSpaceBasis{\lambda}, \qquad \lambda \in \Ptns_{L}. $$
Let $\VertForm\cdot\cdot$ denote the non-degenerate form on $\SymSpace{L}$ defined by the forms (\ref{NondegFormExtn}) on $\RootSpace\g\alpha \otimes \Indet{t}^d$ and by Propositions \ref{TensorFormPropn} and \ref{SymPowerPropn}.  So
$$
\VertForm{\SymSpaceBasis{\lambda}}{\SymSpaceBasis{\mu}} = \prod_{0 \leqslant d \leqslant \Nilp} \prod_{\alpha \in \PositiveRoots} \VertForm{\SymSpaceBasis{\lambda^{\alpha, d}}}{\SymSpaceBasis{\mu^{\alpha, d}}},
$$
for all $\lambda, \mu \in \Ptns_{L}$, where $\VertForm{\SymSpaceBasis{\lambda^{\alpha, d}}}{\SymSpaceBasis{\mu^{\alpha, d}}}$ is defined by Proposition \ref{SymPowerPropn}.
Let $\PtnSpanForm{L}$ be the bilinear form on $\Span(\Ptns_L)$ given by bilinear extension of
$$ \PtnSpanForm{L}\PairForm{\y{\lambda}}{\y{\mu}}{} = \VertForm{\SymSpaceBasis{\lambda}}{\SymSpaceBasis{\mu}}{}, \qquad \lambda, \mu \in \Ptns_L.$$
Then the form $\PtnSpanForm{L}$ is non-degenerate.

\begin{theorem}\label{FormEqualityThm}
For any $L \in \MatrixSpace$, 
$$ \ModShap|_{\Span(\Ptns_L)} = \hElement{L} \cdot \PtnSpanForm{L}$$
where $\hElement{L} \in \SymAlg{\hHat}$ is given by
$$ \hElement{L} = \prod_{0 \leqslant d \leqslant \Nilp} \prod_{\alpha \in \PositiveRoots} (\Factorial{L_{\alpha, d}}) \hspace{0.2em} ( \hElement{\alpha} \otimes \Indet{t}^\Nilp)^{L_{\alpha, d}}.$$
\end{theorem}
\begin{proof}
Let $\lambda, \mu \in \Ptns_L$.  Then:
\begin{eqnarray*}
\Form{\lambda}{\mu} &=& \prod_{0 \leqslant d \leqslant \Nilp} \prod_{\alpha \in \PositiveRoots} \SymPtnPair{\lambda^{\alpha, d}}{\Star{(\mu^{\alpha,d})}} \\
	&& \Note{by Proposition \ref{ExactValuePropn}} \\
	&=& 	\prod_{0 \leqslant d \leqslant \Nilp} \prod_{\alpha \in \PositiveRoots} (\Factorial{L_{\alpha, d}}) \hspace{0.2em} ( \hElement{\alpha} \otimes \Indet{t}^\Nilp)^{L_{\alpha, d}} \VertForm{\SymSpaceBasis{\lambda^{\alpha, d}}}{\SymSpaceBasis{\mu^{\alpha,d}}} \\
	&& \Note{by Proposition \ref{SymPtnPairValuePropn}}\\
	&=& \big[ \prod_{0 \leqslant d \leqslant \Nilp} \prod_{\alpha \in \PositiveRoots} (\Factorial{L_{\alpha, d}}) \hspace{0.2em} ( \hElement{\alpha} \otimes \Indet{t}^\Nilp)^{L_{\alpha, d}} \big]\hspace{0.2em} \cdot \VertForm{\SymSpaceBasis{\lambda}}{\SymSpaceBasis{\mu}}{} \\
	&=& \hElement{L} \cdot \PtnSpanForm{L}\PairForm{\y{\lambda}}{\y{\mu}}{} .
\end{eqnarray*}
The set $\set{\y{\lambda}| \lambda \in \Ptns_L}$ forms a basis for $\Span(\Ptns_L)$, and so the equality follows.
\end{proof}
\end{subsection}

\end{section}

\begin{section}{Reducibility Conditions for Verma Modules}
Let $(\g, \hZero, \h, \g_{+}, \Invol)$ denote a Lie algebra with triangular decomposition and non-degenerate pairing, and let $(\gHat, \hZero, \hHat, \gHat_{+}, \Invol)$ denote the truncated current Lie algebra of nilpotency index $\Nilp$ associated to $\g$.
In this section we establish reducibility criterion for a Verma module $\Verma{\Lambda}$ for $\gHat$ in terms of evaluations of the functional $\Lambda \in \LinearDual\hHat$. 
We then interpret this result separately for the semisimple finite-dimensional Lie algebras, for the affine Kac-Moody Lie algebras, for the symmetrizable Kac-Moody Lie algebras, for the Virasoro algebra and for the Heisenberg algebra.

\begin{theorem}\label{ReducibilityTheorem}
Let $\Lambda \in \LinearDual{\hHat}$ and let $\chi \in \QPlus$.
\begin{enumerate}
\item The Verma module $\Verma{\Lambda}$ for $\gHat$ contains a non-zero primitive vector of weight $\Lambda|_{\hZero} - \chi$ if and only if 
\begin{equation}\label{ReducibilityEqn}
\Functional{\Res\Lambda}{\hElement{\alpha}} = 0
\end{equation}
for some $\alpha \in \PositiveRoots$ such that $\chi - \alpha \in \QPlus$;
\item $\Verma{\Lambda}$ is reducible if and only if equation (\ref{ReducibilityEqn}) holds for some $\alpha \in \Roots$.
\end{enumerate}
\end{theorem}
\begin{proof}
Let $\Lambda \in \LinearDual{\hHat}$ and let $\alpha \in \PositiveRoots$.
By Proposition \ref{RadicalCoincidesWithMaxSubmodule}, the Verma module $\Verma{\Lambda}$ has a non-zero primitive vector of weight $\Lambda|_{\hZero} - \chi$ if and only if the form $\ShapForm_\chi (\Lambda)$ is degenerate.
The determinants $\Det{\ShapForm_\chi}$ and $\Det{\ModShap_\chi}$ can differ only in sign, and 
$$ \Det{\ShapForm_\chi (\Lambda)} = \Functional{\Lambda}{\Det{\ShapForm_\chi}}.$$
Hence such a primitive vector exists if and only if $\Functional{\Lambda}{\Det{\ModShap_\chi}}$ vanishes.  Now
\begin{eqnarray*}
\Functional{\Lambda}{\Det{\ModShap_\chi}}
	&=& \Functional{\Lambda}{\prod_{L \in \MatrixSpace_\chi} \Det{\ModShap|_{\Span(\Ptns_L)}}} \\
	&& \Note{by Corollary \ref{DetBlockCorollary}} \\
	&=& \Functional{\Lambda}{\prod_{L \in \MatrixSpace_\chi} \Det{\PtnSpanForm{L}} \cdot \hElement{L}^{\left| \Ptns_L \right|}} \\
	&& \Note{by Theorem \ref{FormEqualityThm}} \\
	&=& \prod_{L \in \MatrixSpace_\chi} \Det{\PtnSpanForm{L}} \cdot \Functional{\Lambda}{\hElement{L}}^{\left| \Ptns_L \right|}.
\end{eqnarray*}
For any $L \in \MatrixSpace_\chi$, the form $\PtnSpanForm{L}$ is non-degenerate, and so $\Det{\PtnSpanForm{L}}$ is a non-zero scalar.
Hence $\Functional{\Lambda}{\Det{\ModShap_\chi}}$ vanishes if and only if $\Functional{\Lambda}{\hElement{L}}$ vanishes for some $L \in \MatrixSpace_\chi$.
As
$$ \hElement{L} = \prod_{0 \leqslant d \leqslant \Nilp} \prod_{\alpha \in \PositiveRoots} (\Factorial{L_{\alpha, d}}) \hspace{0.2em} ( \hElement{\alpha} \otimes \Indet{t}^\Nilp)^{L_{\alpha, d}},$$
$\Functional{\Lambda}{\Det{\ModShap_\chi}}$ vanishes if and only if $\Functional{\Lambda}{\hElement{\alpha}\otimes \Indet{t}^\Nilp}$ is zero for some $\alpha \in \PositiveRoots$ for which there exists $L \in \MatrixSpace_\chi$ and $0 \leqslant d \leqslant \Nilp$ with $L_{\alpha, d} > 0$.
This condition on $\alpha$ is equivalent to requiring that there exist some partition $\mu \in \Ptns_\chi$ for which $\Length{\mu^\alpha} > 0$, which occurs precisely when $\chi - \alpha \in \QPlus$.
Hence the first part is proven; as $\hElement\alpha$ and $\hElement{-\alpha}$ are proportional, for any $\alpha \in \Roots$, the second part follows.
\end{proof}

It is apparent from Theorem \ref{ReducibilityTheorem} that the reducibility of a Verma module $\Verma{\Lambda}$ for $\gHat$ depends only upon $\Res\Lambda$.

\begin{subsection}{Symmetrizable Kac-Moody Lie Algebras}
Let $\g$ be a symmetrizable Kac-Moody Lie algebra as per Examples \ref{KMLAExample} and \ref{KMLAPairingExample}.
The map
$$ \hElementMap : \LinearDual{\h} \to \h,$$
from Example \ref{KMLAPairingExample} transports the non-degenerate form $\VertForm{\cdot}{\cdot}$ on $\h$ to the space $\LinearDual\h$ via 
$$ \VertForm{\chi}{\gamma} = \VertForm{\hElement{\chi}}{ \hElement{\gamma}}, \qquad \chi, \gamma \in \LinearDual{\h}.$$
Hence, for any $\Lambda \in \LinearDual{\hHat}$ and $\alpha \in \Roots$,
$$
\Functional{\Res\Lambda}{\hElement\alpha} = \VertForm{ \hElement{\Res\Lambda}}{ \hElement\alpha} = \VertForm{\Res\Lambda}{\alpha},
$$
by definition of the map $\hElementMap$.
The following Corollary of Theorem \ref{ReducibilityTheorem} may be viewed as a generalization of Corollary \ref{FDLACorollary}. 
\begin{corollary}\label{KMLACorollary}
Let $\g$ be a symmetrizable Kac-Moody Lie algebra, and let $\gHat$ denote the truncated current Lie algebra of nilpotency index $\Nilp$ associated to $\g$.
Then, for any $\Lambda \in \LinearDual\hHat$, the Verma module $\Verma{\Lambda}$ for $\gHat$ is reducible if and only if $\Res\Lambda$ is orthogonal to some root of $\g$ with respect to the symmetric bilinear form.
\end{corollary}
\end{subsection}

\begin{subsection}{Finite-dimensional semisimple Lie algebras}
The following Corollary is a special case of Corollary \ref{KMLACorollary}.
\begin{corollary}\label{FDLACorollary}
Let $\g$ be a finite-dimensional semisimple Lie algebra, and let $\gHat$ denote the truncated current Lie algebra of nilpotency index $\Nilp$ associated to $\g$.
Then, for any $\Lambda \in \LinearDual\hHat$, the Verma module $\Verma{\Lambda}$ for $\gHat$ is reducible if and only if $\Res\Lambda$ is orthogonal to some root of $\g$ in the geometry defined by the Killing form.
\end{corollary}

Hence the reducibility criterion for Verma modules for $\gHat$ can be described by a finite union of hyperplanes in $\LinearDual\h$.
\begin{example}
Figures \ref{A2Figure} and \ref{G2Figure} illustrate the reducibility criterion for the Lie algebras $\g$ over $\RealNumbers$ with root systems $\FiniteSeries{A}{2}$ and $\FiniteSeries{G}{2}$, respectively.
Roots are drawn as arrows.
A Verma module $\Verma{\Lambda}$ for $\gHat$ is reducible if and only if $\Res\Lambda$ belongs to the union of hyperplanes indicated.
\end{example}

\begin{figure}
\begin{minipage}[t]{8.0cm}
\includegraphics{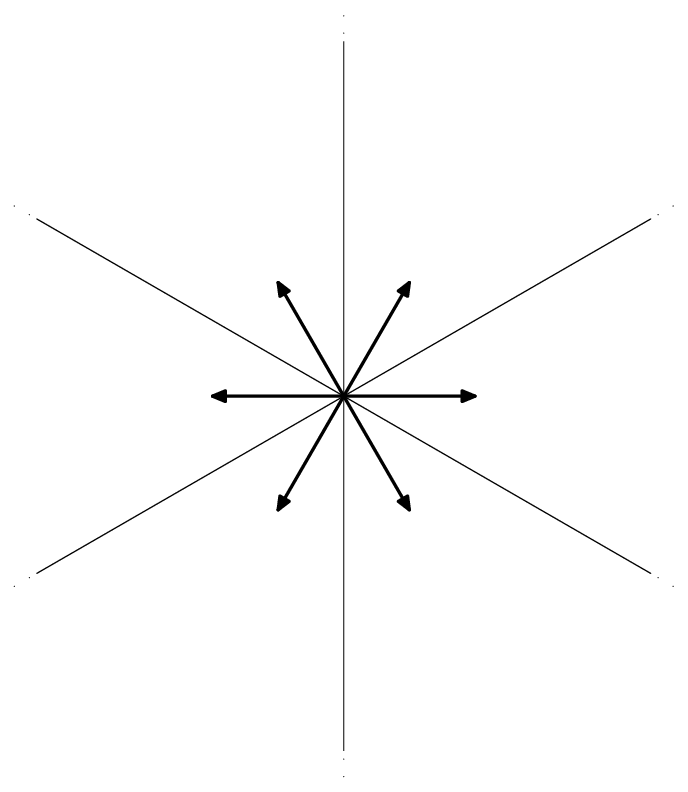}\par
\caption{\label{A2Figure}Criterion for $\FiniteSeries{A}{2}$}
\end{minipage}
\begin{minipage}[t]{8.0cm}
\includegraphics{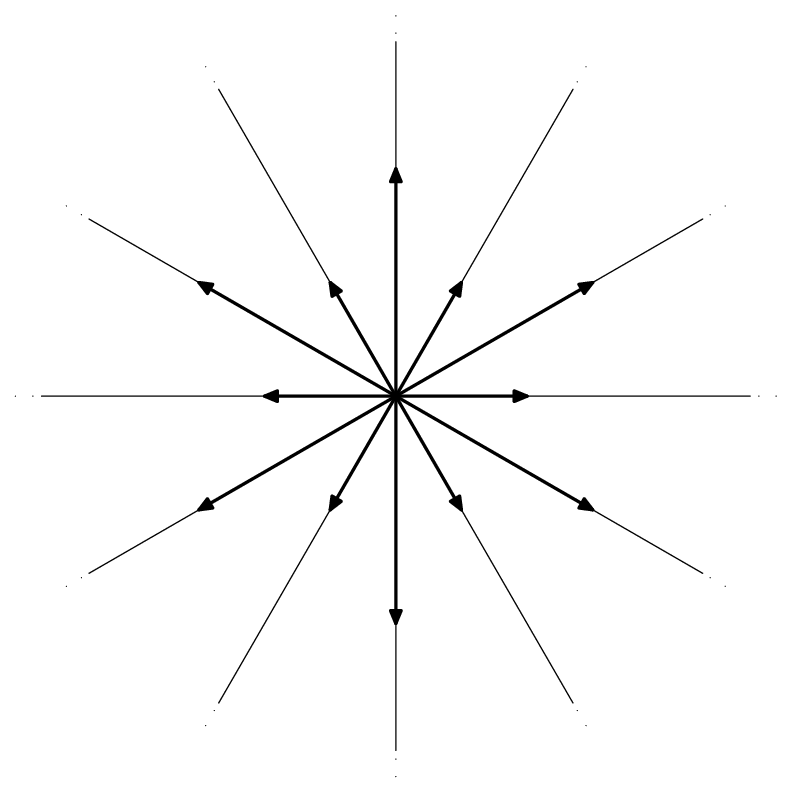}\par
\caption{\label{G2Figure}Criterion for $\FiniteSeries{G}{2}$}
\end{minipage}
\end{figure}
\end{subsection}

\begin{subsection}{Affine Kac-Moody Lie algebras}\label{AffineSubsection}
We refine the criterion of Corollary \ref{KMLACorollary} for the affine Kac-Moody Lie algebras.
Let $\Finite{\g}$ denote a finite-dimensional semisimple Lie algebra over the field $\K$ with Cartan subalgebra $\Finite{\h}$, root system $\Finite\Roots$ and Killing form $\VertForm\cdot\cdot$.
Let $\g$ denote the affinization of $\Finite\g$,
$$ \g = \Finite\g \otimes \K [ \Indet{s}, \Indet{s}^{-1} ] \oplus \K \AffineC \oplus \K \AffineD,$$
with Lie bracket relations
$$
\LieBrac{ x \otimes \Indet{s}^m }{y \otimes \Indet{s}^n} = \LieBrac{x}{y} \otimes \Indet{s}^{m+n} + m \delta_{m, -n} \VertForm{x}{y} \AffineC, 
\qquad \LieBrac{\AffineD}{x \otimes \Indet{s}^m} = m x \otimes \Indet{s}^m,
\qquad \LieBrac{\AffineC}{\g} = 0,
$$
for all $x,y \in \Finite\g$ and $m,n \in \Z$.
Let $\Roots$ denote the root system of $\g$, and let 
$$\h = \Finite\h \oplus \K \AffineC \oplus \K \AffineD,$$
denote the Cartan subalgebra.
Consider any $\Lambda \in \LinearDual{\Finite\h}$ as a functional on $\h$ by declaring
$$ \Lambda (\AffineC) = \Lambda (\AffineD) = 0.$$
This identifies $\LinearDual{\Finite\h}$ with a subspace of $\LinearDual\h$.
Let $\ImagRoot, \OtherRoot \in \LinearDual\h$ be given by
\begin{eqnarray*}
\Functional{\ImagRoot}{\Finite\h} &=& 0, \qquad \Functional\ImagRoot\AffineC = 0, \qquad \Functional\ImagRoot\AffineD = 1, \\
\Functional\OtherRoot{\Finite\h} &=& 0, \qquad \Functional\OtherRoot\AffineC = 1, \qquad \Functional\OtherRoot\AffineD = 0,
\end{eqnarray*}
so that
\begin{equation}\label{DualCartanDecomp}
\LinearDual\h = \LinearDual{\Finite\h} \oplus \K \ImagRoot \oplus \K \OtherRoot.
\end{equation}
The symmetric bilinear form $\VertForm\cdot\cdot$ on $\LinearDual\h$ may be obtained as an extension of the Killing form on $\LinearDual{\Finite\h}$, via
\begin{equation}\label{FormExtension}
\VertForm\ImagRoot{\LinearDual{\Finite\h}} = \VertForm\OtherRoot{\LinearDual{\Finite\h}} = 0, \qquad \VertForm\ImagRoot\ImagRoot = \VertForm\OtherRoot\OtherRoot = 0, \qquad \VertForm\ImagRoot\OtherRoot = 1. 
\end{equation}
The sum (\ref{DualCartanDecomp}) is orthogonal with respect to this form.
For any $\Lambda \in \LinearDual\h$, let $\FiniteRes\Lambda \in \LinearDual{\Finite\h}$ denote the projection of $\Lambda$ on to $\LinearDual{\Finite\h}$ defined by the decomposition (\ref{DualCartanDecomp}).
The root system $\Roots = \RealRoots \cup \ImaginaryRoots$ of $\g$ is given by, 
\begin{equation}\label{AffineRoots}
\RealRoots = \set{ \alpha + m \ImagRoot | \alpha \in \Finite\Roots, \ m \in \Z}, \qquad \ImaginaryRoots = \set{ m \ImagRoot | m \in \Z, \ m \ne 0}.
\end{equation}
\begin{corollary}\label{AffineCorollary}
Let $\g$ denote an affine Kac-Moody Lie algebra, and let $\gHat$ denote the truncated current Lie algebra of nilpotency index $\Nilp$ associated to $\g$.
Then, for any $\Lambda \in \LinearDual\hHat$, the Verma module $\Verma\Lambda$ for $\gHat$ is reducible if and only if $\Functional{\Res\Lambda}\AffineC = 0$ or
$ \VertForm{\FiniteRes{\Res\Lambda}}{\alpha} = m \hspace{0.15em} \Functional{\Res\Lambda}\AffineC$
for some $\alpha \in \Finite\Roots$ and $m \in \Z$.
\end{corollary}
\begin{proof}
It is immediate from (\ref{DualCartanDecomp}) and (\ref{FormExtension}) that
$$ \Lambda = \FiniteRes\Lambda + \VertForm\Lambda\OtherRoot \ImagRoot + \VertForm\Lambda\ImagRoot \OtherRoot.$$
Hence
$$
\Functional\Lambda\AffineC = \VertForm\Lambda\ImagRoot \Functional\OtherRoot\AffineC = \VertForm\Lambda\ImagRoot.$$
Therefore $\Functional\Lambda\AffineC = 0$ if and only if $\VertForm\Lambda\beta = 0$ for some $\beta \in \ImaginaryRoots$.
For $\alpha \in \Finite\Roots$ and $m \in \Z$,
$$
\VertForm\Lambda{\alpha + m \ImagRoot} = \VertForm\Lambda\alpha + m \VertForm\Lambda\ImagRoot = \VertForm{\FiniteRes\Lambda}{\alpha} + m \Functional\Lambda\AffineC
$$
and so $\VertForm\Lambda{\alpha + m \ImagRoot} = 0$ if and only if
$ \VertForm{\FiniteRes\Lambda}{\alpha} = -m \Functional\Lambda\AffineC$.
The claim now follows from (\ref{AffineRoots}) and Corollary \ref{KMLACorollary}.
\end{proof}
\begin{example}
Figures \ref{AffineA2Figure}, \ref{AffineB2Figure} and \ref{AffineG2Figure} illustrate the reducibility criterion of Corollary \ref{AffineCorollary} in the case where $\K = \RealNumbers$ and $\Finite\g$ is the Lie algebra with root systems $\FiniteSeries{A}{2}$, $\FiniteSeries{B}{2}$ and $\FiniteSeries{G}{2}$, respectively.
A Verma module $\Verma{\Lambda}$ for $\gHat$ is reducible if and only if ${\FiniteRes{\Res\Lambda}}$ belongs to the described infinite union of hyperplanes, where the length of the dashed line segment is $| \Functional{\Res\Lambda}{\AffineC} |$ times the length of a short root for $\Finite\g$.

\begin{figure}
\begin{minipage}[t]{8.0cm}
\includegraphics{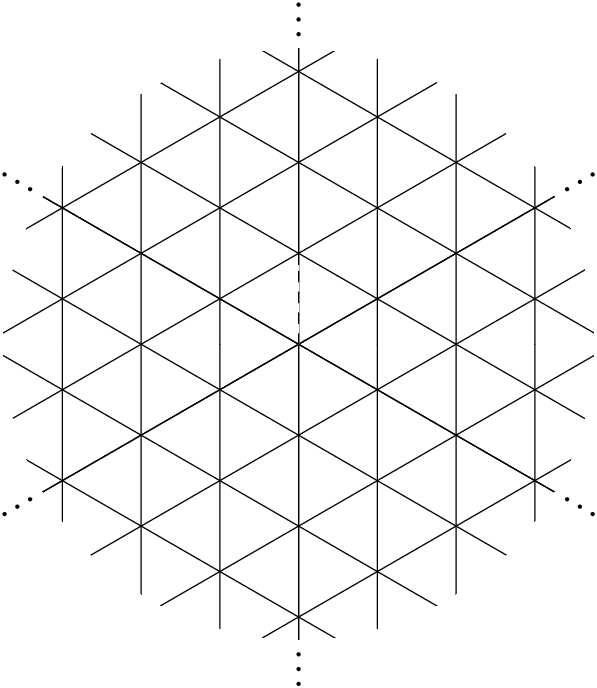}\par
\caption{\label{AffineA2Figure}Criterion for affine $\FiniteSeries{A}{2}$}
\end{minipage}
\begin{minipage}[t]{8.0cm}
\includegraphics{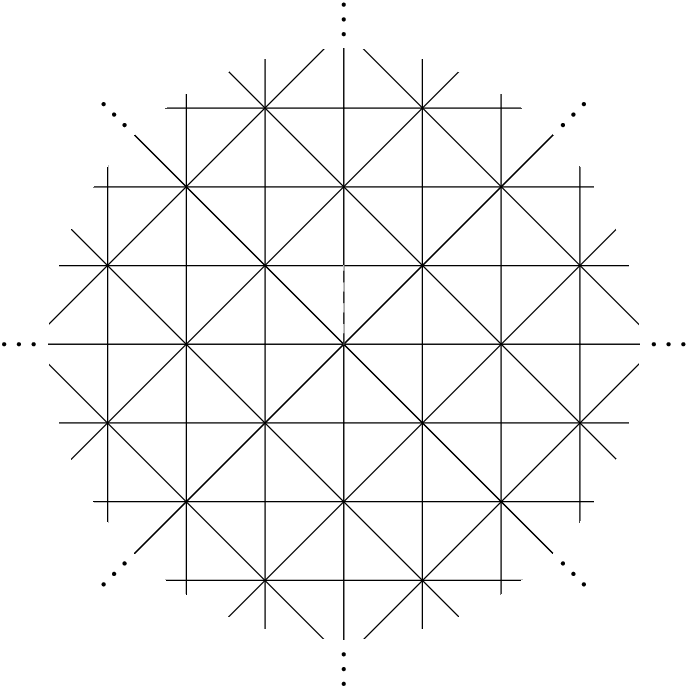}\par
\caption{\label{AffineB2Figure}Criterion for affine $\FiniteSeries{B}{2}$}
\end{minipage}
\end{figure}
\begin{figure}
\includegraphics{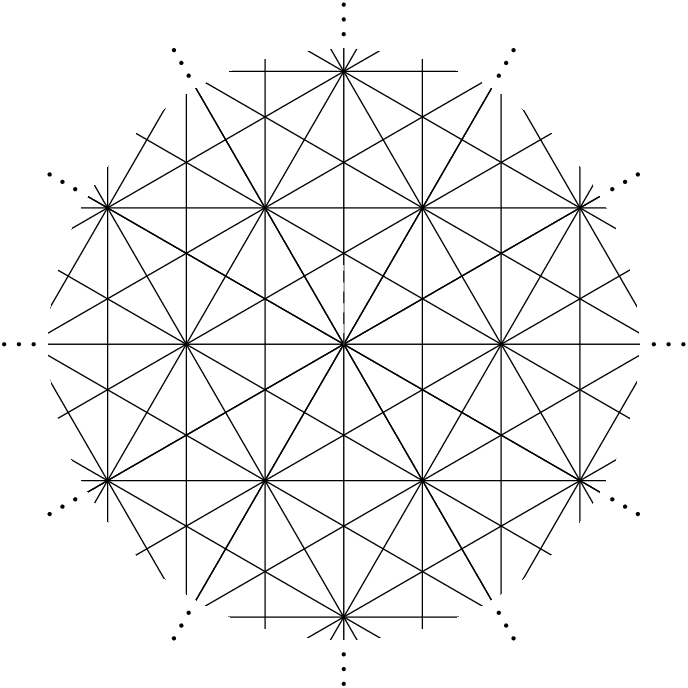}\par
\caption{\label{AffineG2Figure}Criterion for affine $\FiniteSeries{G}{2}$}
\end{figure}
\end{example}
\end{subsection}
\begin{subsection}{The Virasoro Algebra}
The following Corollary is immediate from Theorem \ref{ReducibilityTheorem} and Examples \ref{VirasoroExample} and \ref{VirasoroPairingExample}.
\begin{corollary}\label{VirasoroCorollary}
Let $\g$ denote the Virasoro algebra, and let $\gHat$ denote the truncated current Lie algebra of nilpotency index $\Nilp$ associated to $\g$.
Then, for any $\Lambda \in \LinearDual\hHat$, the Verma module $\Verma{\Lambda}$ for $\gHat$ is reducible if and only if 
$$ 2m \Functional{\Res\Lambda}{\VirL_0} + \psi(m) \Functional{\Res\Lambda}{\VirC} = 0,$$
for some non-zero integer $m$.
\end{corollary}
Hence, if $\psi$ is defined by $\psi(m) = \SFrac{m^3 -m}{12}$ and $\K= \RealNumbers$, a Verma module $\Verma{\Lambda}$ for $\gHat$ is reducible if and only if $\Res\Lambda$ belongs to the infinite union of hyperplanes indicated in Figure \ref{VirasoroFigure}.
The extension of a functional in the horizontal and vertical directions is determined by evaluations at $\VirC$ and $\VirL_0$, respectively.

\begin{figure}
\includegraphics{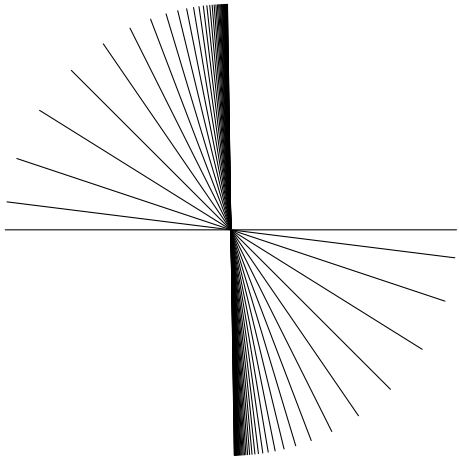}\par
\caption{\label{VirasoroFigure}Criterion for the Virasoro algebra}
\end{figure}
\end{subsection}
\begin{subsection}{The Heisenberg Algebra}
The following Corollary is immediate from Theorem \ref{ReducibilityTheorem} and Examples \ref{HeisenbergExampleOne} and \ref{HeisenbergPairingExample}.
\begin{corollary}\label{HeisenbergCorollary}
Let $\hat{\Heis}$ denote the truncated current Lie algebra of nilpotency index $\Nilp$ associated to the Heisenberg algebra $\Heis$.
Then, for any $\Lambda \in \LinearDual{\hHat}$, a Verma module $\Verma\Lambda$ for $\hat{\Heis}$ is reducible if and only if $\Functional{\Res\Lambda}{\HeisH} = 0$.
\end{corollary}
\end{subsection}
\end{section}

\begin{section}{Acknowledgements}
The author is most grateful for the patient guidance and encouragement of Yuly Billig during his recent visit to the University of Sydney.
This work was completed as part of a cotutelle Ph.~D.~ programme under the supervision of Vyacheslav Futorny at the Universidade de S\~ao Paulo and Alexander Molev at the University of Sydney.
\end{section}

\appendix
\begin{section}{Characters of Irreducible Highest-Weight Modules for Truncated Current Lie Algebras}
Let $\g$ denote a Lie algebra with triangular decomposition and non-degenerate pairing,
and let $\gHat$ denote the truncated current Lie algebra of nilpotency index $\Nilp$ associated to $\g$.
Theorem \ref{ReducibilityTheorem} describes a reducibility criterion for Verma modules for $\gHat$, but provides little information on the size of the maximal submodule.
This appendix describes the characters of the irreducible highest-weight $\gHat$-modules under the assumption that $\h$ is one-dimensional (and hence $\hZero = \h$).
For example, $\g$ may be the Lie algebra $\SpecialLinear{2}$, the Witt algebra, or a modified Heisenberg algebra.

For any $\gamma \in \LinearDual\h$, let $\Irred\gamma$ denote the irreducible highest-weight $\g$-module of highest-weight $\gamma$, and
for any $\Lambda \in \LinearDual{\hHat}$, let $\Irred\Lambda$ denote the irreducible highest-weight $\gHat$-module of highest weight $\Lambda$.
Let 
$$ \set{\Exp\chi | \chi \in \LinearDual{\h}}$$
denote a multiplicative copy of the additive group $\LinearDual\h$, so that
$$ \Exp\chi \cdot \Exp\gamma = \Exp{\chi + \gamma}, \qquad \chi, \gamma \in \LinearDual\h.$$
If $M$ is a vector space graded by $\LinearDual\h$,
$ M = \oplus_{\chi \in \LinearDual\h} \WeightSpace{M}{\chi}$,
such that all components $\WeightSpace{M}{\chi}$ are finite-dimensional, write
$$ \Character{M} = \sum_{\chi \in \LinearDual\h} (\Dim \WeightSpace{M}\chi) \hspace{0.15em} \Exp\chi.$$
\begin{propn}
Let $\g, \gHat$ be as above, and let $\Lambda \in \LinearDual{\hHat}$.
Let $0 \leqslant m \leqslant \Nilp$ be minimal such that $\Lambda_n = 0$ for all $m < n \leqslant \Nilp$.
Then, if $m > 0$,
$$ \Character{\Irred\Lambda} = \Exp{\Lambda_0} \cdot (\Character{\UEA{\g_{-}}})^m,$$
and if $m = 0$, then $\Irred\Lambda$ is a $\g$-module isomorphic to $\Irred{\Lambda_0}$.
\end{propn}
\begin{proof}
Suppose that $m=0$.
Since $\g$ is the quotient of $\gHat$ by the ideal $\oplus_{0<i \leqslant \Nilp} \g \otimes \Indet{t}^i$, the $\g$-module $\Irred{\Lambda_0}$ is a natural $\gHat$-module.
Moreover, $\Irred{\Lambda_0}$ is an irreducible highest-weight $\gHat$-module of highest-weight $\Lambda$, and so $\Irred{\Lambda_0} \cong \Irred{\Lambda}$.

Suppose instead that $m > 0$.
Then it must be that $\Lambda_m \ne 0$.
Let $\gHat'$ denote the truncated current Lie algebra of nilpotency index $m$ associated to $\g$.  Let
$$ \Lambda' = (\Lambda_0, \dots, \Lambda_m) \quad \in \LinearDual{(\hHat')},$$
and let $\Verma{\Lambda'}$ denote the Verma module for $\gHat'$ of highest-weight $\Lambda'$.
Since $\gHat'$ is the quotient of $\gHat$ by the ideal $\oplus_{m < i \leqslant \Nilp} \g \otimes \Indet{t}^i$, the $\gHat'$-module $\Verma{\Lambda'}$ is a natural $\gHat$-module.
Moreover, $\Verma{\Lambda'}$ is of highest-weight $\Lambda$ as a $\gHat$-module.
Since $\h$ is one-dimensional,  $\Verma{\Lambda'}$ is an irreducible $\gHat'$-module, by Theorem \ref{ReducibilityTheorem}.
Hence $\Verma{\Lambda'}$ is the irreducible $\gHat$-module of highest-weight $\Lambda$, \ie $\Irred{\Lambda} \cong \Verma{\Lambda'}$ as $\gHat$-modules.
In particular, $\Irred\Lambda$ and $\Verma{\Lambda'}$ are isomorphic as $\LinearDual{\hHat}$-graded vector spaces.  Now
$$ \Character{\Verma{\Lambda'}} = \Exp{\Lambda_0} \cdot (\Character{\UEA{\g_{-}}})^m$$
by Proposition \ref{VermaPropertiesPropn} part \ref{VermaPropertiesPropn2}, and so the claim follows.
\end{proof}

\end{section}

\begin{section}{Imaginary Highest-Weight Theory for Truncated Current Lie Algebras}
Let $\Finite\g$ denote the finite-dimensional Lie algebra $\SpecialLinear{2}$ over the field $\K$, with root system $\Finite\Roots = \set { \pm \upalpha}$, and let $\g$ denote the affinization of $\Finite\g$ (\cf Subsection \ref{AffineSubsection}).
Let $\h$ denote the Cartan subalgebra of $\g$, let $\Roots$ denote the root system, and let $\ImagRoot$ denote the fundamental imaginary root.  
Let
\begin{equation*}
\PositiveRoots = \set{ \upalpha + m \ImagRoot | m \in \Z} \cup \set{ m \ImagRoot | m \in \Z, \ m > 0},
\end{equation*}
so that $\Roots = \PositiveRoots \cup - \PositiveRoots$.
Let
$$ \g_{+} = \oplus_{\beta \in \PositiveRoots} \RootSpace\g\beta, \qquad \g_{+} = \oplus_{\beta \in \PositiveRoots} \RootSpace\g{-\beta},$$
so that
\begin{equation}\label{ImaginaryDecomposition}
\g = \g_{-} \oplus \h \oplus \g_{+}.
\end{equation}
The subset $\PositiveRoots \subset \Roots$ is called the \NewTerm{imaginary partition} of the root system.
The decomposition (\ref{ImaginaryDecomposition}) defined by $\PositiveRoots$ does not satisfy the axioms of a triangular decomposition in the sense of Section \ref{LATDSection}, nor in the sense of \cite{MoodyPianzola}: 
the additive semigroup $\QPlus$ generated by $\PositiveRoots$ is not generated by any linearly independent subset of $\QPlus$.
Nevertheless, there exists a non-classical highest-weight theory for $\g$ defined by the decomposition (\ref{ImaginaryDecomposition}).
This \NewTerm{imaginary highest-weight theory}, pioneered by Futorny, is remarkably different from classical highest-weight theory (\cite{FutornyAffineOne}, \cite{FutornyAffineTwo}, \cite{WilsonImaginaryHWT}).
In particular, weight spaces of imaginary highest-weight modules may be infinite-dimensional.

Let $\gHat$ denote the truncated current Lie algebra of nilpotency index $\Nilp$ associated to $\g$.
We investigate the difficulty inherent in employing our techniques to derive reducibility criterion for the Verma modules $\Verma\Lambda$ for $\gHat$, $\Lambda \in \LinearDual\hHat$.
The Theorem \ref{FormEqualityThm} holds in this setting.
However, as we shall see, the degeneracy of an evaluation 
$ \Functional\Lambda{\ModShap_\chi}$ of the modified Shapovalov form 
may not be deduced from the degeneracy of the evaluations $\Functional\Lambda{\ModShap|_{\Span(\Ptns_L)}}$, where $L \in \Ptns_\chi$.

Let $\Nilp = 1$, and as per Sections \ref{HWTSection} and \ref{TCLASection}, let
$$\C = \PositiveRoots, \qquad \CHat = \C \times \set{0,1}.$$
All root spaces of $\g$ are one-dimensional, so the choice of basis for $\g_{+}$
$$ \C \ni \quad \beta \quad \leftrightarrow \quad \x\beta \quad \in \RootSpace\g\beta$$
is unique up to scalar multiples.
Fix the order of the basis elements $\set{\x\gamma | \gamma \in \CHat}$ by firstly comparing degree in the indeterminate $\Indet{t}$, and secondly by the following order of $\set{\x\beta | \beta \in \C}$:
$$ \cdots \x{\upalpha - 2 \ImagRoot}, \ \x{\upalpha - \ImagRoot}, \ \x{\upalpha}, \ \x{\upalpha + \ImagRoot}, \ \x{\upalpha + 2 \ImagRoot}, \cdots 
\qquad \cdots \x\ImagRoot, \ \x{2 \ImagRoot}, \cdots$$
For any integer $m>0$, define partitions
$$ \mu_m = \set{ (\upalpha - m \ImagRoot, 0)} \cup \set{ (\ImagRoot, 0) \ \ \text{($m$ times)}}, \qquad \gamma_m = \set{ (\upalpha - m \ImagRoot, 1)} \cup \set{ (\ImagRoot, 1) \ \ \text{($m$ times)}},$$
and $\lambda = \set{ (\upalpha, 0)}$.
Let $\chi = \upalpha \in \QPlus$.  Then
$$ \set{\lambda} \cup \set{ \mu_m, \gamma_m | m > 0} \subset \Ptns_\chi.$$
Elementary computation using the Lie bracket relations shows that, for any $m >0$, 
$$ \Form{\lambda}{\gamma_m} = (-2)^m \hElement{\upalpha - m \ImagRoot} \otimes \Indet{t}^0, \qquad \Form{\mu_m}{\lambda} = (-2)^m \hElement{\upalpha -m \ImagRoot} \otimes \Indet{t}^1.$$
Hence, if the basis $\set{ \y\mu | \mu \in \Ptns_\chi}$ is to be linearly ordered so that the matrix representation of $\ModShap_\chi$ is upper triangular, then it must be that both
$$ \y{\mu_m} < \y\lambda, \qquad \y\lambda < \y{\gamma_m},$$
for all $m>0$.
Thus the matrix of $\ModShap_\chi$ would be bilaterally infinite.

The degeneracy of a bilaterally-infinite upper-triangular matrix can not be determined from its diagonal entries, as the following simple example demonstrates.
Let $V$ denote the vector space with basis the symbols 
\begin{equation}\label{VBasis}
\set{ \mathrm{v_m} | m \in \Z}
\end{equation}
and let $\Phi : V \to V$ be defined by linear extension of the rule
$$ \Phi : \mathrm{v}_m \mapsto \mathrm{v}_{m+1}, \qquad m \in \Z.$$
Then $\Phi$ is an automorphism of $V$.
Order the basis elements (\ref{VBasis}) by their indices.
Then the matrix representation $M$ of $\Phi$ with respect to this ordered basis will be upper triangular, in the sense that $M_{i,j} = 0$ whenever $i>j$.
However, the diagonal entries $M_{i,i}$ are all identically zero.
\end{section}

\bibliographystyle{abbrv}
\bibliography{thesis}
\end{document}